\documentclass[12pt,leqno]{amsart}

\usepackage{graphicx}
\usepackage{amssymb}
\newtheorem{thm}{Theorem}
\newtheorem{cor}[thm]{Corollary}
\newtheorem{lem}[thm]{Lemma}
\newtheorem{prop}[thm]{Proposition}
\theoremstyle{definition}
\newtheorem{defn}[thm]{Definition}
\theoremstyle{remark}

\numberwithin{equation}{section}

\newcommand{\To}{\longrightarrow}

\begin{document}
\setcounter{tocdepth}{1}


\title[]{Multiple gaps}
\author{Antonio Avil\'es and Stevo Todorcevic}
\address{Universidad de Murcia, Departamento de Matemáticas, Campus de Espinardo 30100 Murcia, Spain.} \email{avileslo@um.es}
\address{CNRS FRE 3233, Université Paris Diderot Paris 7, UFR de mathématiques case 7012, site Chevaleret
75205 Paris, France. Department of Mathematics, University of Toronto, Toronto, Canada, M5S 3G3.}%
\email{stevo@logique.jussieu.fr, stevo@math.toronto.edu}
\thanks{A. Avil\'{e}s was supported by MEC and FEDER (Project MTM2008-05396), Fundaci\'{o}n S\'{e}neca
(Project 08848/PI/08), Ram\'{o}n y Cajal contract (RYC-2008-02051) and an FP7-PEOPLE-ERG-2008 action. S. Todorcevic was supported by grants from NSERC and CNRS}

\begin{abstract}
We study a higher-dimensional version of the standard notion of a gap formed by a finite sequence of ideals of the quotient algebra $\mathcal{P}(\omega)/fin$. We examine  different types of such objects found in $\mathcal{P}(\omega)/fin$ both from the combinatorial and the descriptive set-theoretic side.
\end{abstract}


\maketitle

\section{Introduction}

Gaps in the quotient algebra $\mathcal{P}(\omega)/fin$ is a phenomenon discovered by Hausdorff more than a century ago.
Their study has always been a prominent theme in set theory and its applications. For example, Kunen's study of gaps of  $\mathcal{P}(\omega)/fin$ 
in models of Martin's axiom has led Woodin to a proof that the Proper Forcing Axiom implies Kaplanski's conjecture about automatic continuity of norms in certain Banach algebras (see, \cite{DW}). Gaps in $\mathcal{P}(\omega)/fin$ were also the main motivation behind the introduction of the set theoretic principle about the chromatic number of open graphs on separable metric spaces (see \cite{TodorcevicTopology}), a principle with many applications (see, for example, \cite{Farah1}, \cite{Farah2}). In \cite{Todorcevicgap}, the second author has initiated the study of gaps in
$\mathcal{P}(\omega)/fin$ from the descriptive set-theoretic side, a study also of independent interest and important applications (see, for example, \cite{TodorcevicRamsey}). The purpose of this paper is to extend this theory to higher dimensions, or in other words, to build a theory that involve gaps formed by more than two ideals of  $\mathcal{P}(\omega)/fin.$ While this theory is not something that immediately suggests itself when one mentions gaps formed by more than two ideals (see, for example, \cite{Farah3}, \cite{Talayco}) it is nevertheless quite natural and it could have been discovered long ago. We came to it while trying to understand the reasons behind the fact that the Banach space $\ell_\infty/c_0 $ is not injective, a result originally due to Amir \cite{Amir}.
While the characterization of $1$-injective spaces due to Kelley \cite{Kelley} suggests that gaps of $\mathcal{P}(\omega)/fin$ should play a role in Amir's result, we were surprised after we realized that the classical theory of gaps in $\mathcal{P}(\omega)/fin$ is actually not useful in its proof.
What is relevant in this context is a notion of a $n$-gap formed by a sequence $\mathcal{I}_i$ $(i<n)$ of ideals of $\mathcal{P}(\omega)/fin$ so that the corresponding sequence of open subsets of $\omega^*=\beta{\omega}\setminus \omega$ has a common boundary point. We give a precise explanation of this in Section 8 of the paper where we use Ditors's analysis of lower bound of norm of averaging operators. This also suggests various refinements of the notion of an $n$-gap that we study below. As said above we study the new notion both from the combinatorial and descriptive set theoretic side.
For example, we show that there is no Hausdorff phenomenon for $n$-gaps where $n>2,$ or more precisely, unlike in the case $n=2,$ the usual axioms of set theory are insufficient for constructing $n$-gaps consisting of $\aleph_1$-generated ideals for $n>2.$ On the other hand, we show that the descriptive set theory of $2$-gaps initiated in \cite{Todorcevicgap} extends naturally to all higher dimensions.

\section{Notation}

For a Boolean algebra $\mathcal{B}$, we denote as $\vee$, $\wedge$ and $\leq$ the join operation, meet operation and order in an abstract Boolean algebra $\mathcal{B}$, and by 0 and 1 its maximum and minimum. Two subsets $I$ and $J$ of $\mathcal{B}$ are called orthogonal if $a\wedge b = 0$ for all $a\in I$, $b\in J$. Other notation that we shall follow are:
$I\vee J = \{a\vee b : a\in I, b\in J\}$ and similarly $I\wedge J$; $I^\perp = \{b\in\mathcal{B} : \forall a\in I\ a\wedge b = 0\}$; $I|_a = \{b\in I : b\leq a\}$; $c\geq I$ if $c\geq a$ for all $a\in I$. Also, $n=\{0,1,\ldots,n-1\}$.\\

Although we give definitions for arbitrary Boolean algebras, which is the natural context, most of the paper deals only with the algebra $\mathcal{P}(\omega)/fin$ and the results refer to this algebra unless otherwise explicitly stated. The elements of $\mathcal{P}(\omega)/fin$ are equivalence classes of the power set of natural numbers $\mathcal{P}(\omega)$ under the relacion of equality modulo finite sets: $a=^\ast b$ if $(a\setminus b)\cup (b\setminus a)$ is finite. However, manipulating equivalence classes is unconvenient for it makes considerable noise in notation. Instead, we prefer to work in $\mathcal{P}(\omega)$ while referring to properties of $\mathcal{P}(\omega)/fin$. So we use the following conventions along the text:\\

Ideals of $\mathcal{P}(\omega)/fin$ are identified with ideals of $\mathcal{P}(\omega)$ that contain the ideal $fin$ of finite sets. Hence, if $I$ is such an ideal and we write $a\in I$, we understand that $a\in P(\omega)$ and $a$ belongs to the corresponding ideal of $P(\omega)$ that contains $fin$. On the other hand, when we state properties of such ideals like being orthogonal, being a multiple gap, a jigsaw, clover, etc. those properties always refer to properties as ideals in the Boolean algebra $\mathcal{P}(\omega)/fin$, never as ideals in $\mathcal{P}(\omega)$. One exception is when we say that an ideal $I$ is analytic or Borel, then we mean as a subset of $\mathcal{P}(\omega) = 2^\omega$.\\

For $a,b\in \mathcal{P}(\omega)$ when we write $a\subset b$ we mean actual inclusion in  $\mathcal{P}(\omega)$ while $a\subset^\ast b$ means that $b\setminus a$ is finite. The use of symbols that refer to the structure of a general Boolean algebra, like  $\leq$ and $\perp$, always refer to $\mathcal{P}(\omega)/fin$, never to $\mathcal{P}(\omega)$. For instance if $a,b\in\mathcal{P}(\omega)$ and $I$ is an ideal of $\mathcal{P}(\omega)/fin$, then $a\leq b$ means $a\subset^\ast b$, $a\perp b$ means $a\cap b =^\ast \emptyset$, $I\leq a$ means $x\subset^\ast a$ for all $x\in I$, $a\in I^\perp$ means that $a\cap x=^\ast\emptyset$ for all $x\in I$, etc. We denote $\omega^\ast = \beta\omega\setminus\omega$ the dual compact space of $\mathcal{P}(\omega)/fin$ under Stone's duality.\\

Given a set $S$, by $S^{<\omega}$ we denote the set of all finite sequences of elements of $S$. Given $s=(s_i)_{i<n}$ and $t=(t_i)_{i<m}$ elements of $S^{<\omega}$, we write $s\leq t$ when $n\leq m$ and $s_i = t_i$ for all $i<n$. We denote $s^\frown t = (s_0,\ldots,s_{n-1},t_0,\ldots,t_{m-1})$. The set $S^{<\omega}$ is a tree when endowed with the partial order $(\leq)$. A tree is a partially ordered set $(T\leq)$ with a minimum element (its root) such that for every $t\in T$ the set $\{s\in T :s<t\}$ is well ordered. A linearly ordered subset of $T$ is called a chain. A branch of $T$ is a maximal chain. An antichain is a subset of $T$ where no two elements satisfy $s<t$. The set of branches of $T$ is denoted by $[T]$. Given $t\in T$, we write $T_t = \{ s\in T : t\leq s\}$ and $[T_t] = \{x\in [T] : t\in x\}$. When $T=S^{<\omega}$ we can identify $[T] = S^\omega$.

Often along the paper, instead of working in $\mathcal{P}(\omega)/fin$ we work in $\mathcal{P}(E)/fin$
where $E$ is a countable set with some particular structure, like $E=T$ being a countable tree or $E=D$ being a countable subset of a compact metric space. The same conventions as for $\mathcal{P}(\omega)/fin$ apply.\\

\section{Basic Notions and Preliminaries}

A gap in a Boolean algebra $\mathcal{B}$ is a couple of orthogonal ideals that cannot be separated by two disjoint elements of $\mathcal{B}$. We introduce a multidimensional generalization:

\begin{defn}
We say that a finite family $\{I_i : i\in n\}$ of mutually orthogonal ideals of a Boolean algebra $\mathcal B$ constitutes a multiple gap, or $n$-gap, if for every function $c:n\To\mathcal{B}$ such that $c(i)\geq I_i$ for all $i\in n$, we have that $\bigwedge_{i\in n}c(i) \neq 0$.\\
\end{defn}

There is a simple topological interpretation of this definition. Let $St(\mathcal B)$ be the Stone dual compact space of the Boolean algebra $\mathcal B$, that is, the space of all ultrafilters of $B$ with the topology generated by the sets $\{\mathfrak{U} : a\in\mathfrak{U}\}$ for $a\in \mathcal B$. For an ideal $I$, let $U(I) = \{\mathfrak{U} : \mathfrak U\cap I \neq\emptyset\}$ be the associated open subset of $St(\mathcal{B})$. The fact that the ideals $\{I_i : i\in n\}$ form an $n$-gap is equivalent to the fact that $\bigcap_{i\in n}\overline{U(I)}\neq\emptyset$. Thus, a multiple gap is nothing else than a finite family of pairwise disjoint open sets whose closures have nonempty intersection.\\

The multiple gap is called dense if $(\bigvee_{i\in n}I_i)^\perp = 0$. Topologically, this means that the union of the open sets considered is dense in $St(\mathcal{B})$. When $n>2$, we can distinguish a variety of types of $n$-gaps depending on how subgaps of lower dimension interact. We point out two extreme types that we found of special interest, and we called clovers and jigsaws.\\

\begin{defn}
Let $\{I_i : i\in n\}$ be a multiple gap and $B$ a nonempty proper subset of $n$. We call the multiple gap a $B$-clover if there exists no $b\in\mathcal{B}$ such that $\{I_j|_b : j\in B\}$ is a multiple gap while $b\in I_i^\perp$ for all $i\not\in B$.
\end{defn}

The multiple gap is called a clover if it is a $B$-clover for all nonempty proper subsets $B$ of $n$. In topological terms, being a $B$-clover is equivalent to the equality $\bigcap_{i\in B}\overline{U(I_i)} \setminus \bigcup_{i\not\in B}\overline{U(I_i)}=\emptyset$. Hence, for instance the picture of a triple clover would look like this:

\begin{center}\includegraphics[scale=0.75]{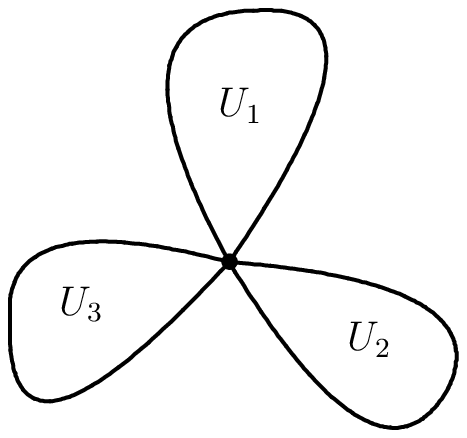}\end{center}

where the intersection of any two closures equals the intersection of all closures (that may consist of more than one point, in spite of the picture).\\

The second type, opposite to clovers is that of jigsaws.\\

\begin{defn}
Let $\{I_i : i\in n\}$ be an $n$-gap, and let $B$ be a nonempty proper subset of $n$. We call the multiple gap a $B$-jigsaw if for every $A\supset B$ and every $a\in\mathcal{B}$, if $\{I_i|_a : i\in A\}$ is a multiple gap, then there exists $b\subset a$ such that $\{I_i|_b : i\in B\}$ is a multiple gap while $b\in I_i^\perp$ for all $i\in A\setminus B$.
\end{defn}

The multiple gap is called a jigsaw if it is a $B$-jigsaw for all nonempty proper subsets $B$ of $n$. In topological terms, being a $B$-jigsaw is equivalent to state that the set $\bigcap_{i\in B}\overline{U(I_i)}\setminus\bigcup_{i\not\in B}\overline{U(I_i)}$ is dense in $\bigcap_{i\in B}\overline{U(I_i)}$. The picture of a triple jigsaw is then the following:

\begin{center}\includegraphics[scale=0.75]{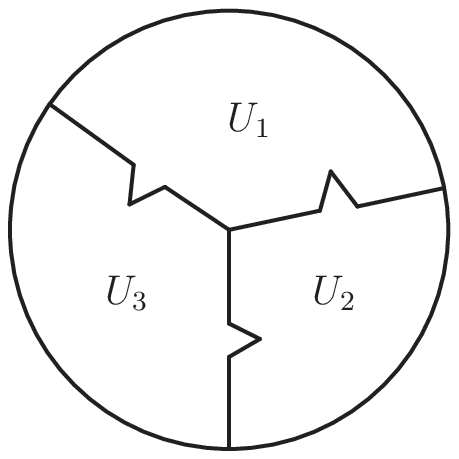}\end{center}

where every point in the intersection of the three closures (in the picture, the central point) can be approached by points that lie in the intersection of exactly two closures.\\

The concept of countable separation plays a central role in our study. Countably separated gaps are introduced in \cite{Todorcevicgap}: Two orthogonal ideals $I_0$ and $I_1$ are countably separated if there exist a sequence of elements $\{c_n : n\in\omega\}$ in the Boolean algebra such that for every $x\in I_0$ and $y\in I_1$ there exists $n$ with $x\leq c_n$ and $y\perp c_n$. Regarding this notion, the following dichotomy is proven in \cite{Todorcevicgap}:

\begin{thm}\label{2gapdichotomy} Let $I_0$ and $I_1$ be two orthogonal analytic ideals in $\mathcal{P}(\omega)/fin$. Then
\begin{enumerate}
\item either they are countably separated in $\mathcal{P}(\omega)/fin$
\item or there exist $\{a_s : s\in 2^\omega\}\subset I_0$ and $\{b_s: s\in 2^\omega\}\subset I_1$ such that $a_s\cap b_s = \emptyset$, $(a_s\cap b_t)\cap (a_t\cap b_s)\neq\emptyset$ for $t\neq s$, and the maps $s\mapsto a_s$ and $s\mapsto b_s$ are continuous.
\end{enumerate}
\end{thm}

Countable separation can be generalized in several ways to multiple gaps, namely:

\begin{defn}
Let $\mathfrak{X}$ be a family of subsets of $n$, and let $\{I_i : i\in n\}$ be ideals in a Boolean algebra. We say that a multiple gap $\{I_i : i<n\}$ is $\mathfrak{X}$-countably separated if there exist elements $\{c_i^k, i\in n, k\in \omega\}$ in the Boolean algebra such that \begin{enumerate} \item $\bigwedge_{i\in A} c_i^k = 0$ for every $k\in\omega$ and every $A\in\mathfrak{X}$, and\item for every $x_0\in I_0$, $\ldots$, $x_{n-1}\in I_{n-1}$, there exists $k\in\omega$ such that $x_i\leq c_i^k$ for all $i$. \end{enumerate}
\end{defn}

The two extreme cases of this definition are of special interest for us (we use the notation $[x]^k = \{z\subset x: |x| = k\}$):

\begin{itemize}
\item We say that the ideals $\{I_i : i\in n\}$ are weakly countably separated if they are $\mathfrak{X}$-countably separated for $\mathfrak{X} = [n]^n$.
\item  We say that the ideals $\{I_i : i\in n\}$ are strongly countably separated if they are $\mathfrak{X}$-countably separated for $\mathfrak{X} = [n]^2$.\\

\end{itemize}

In Section \ref{dichotomy} we prove a generalization of Theorem~\ref{2gapdichotomy} for multiple gaps. We define for every $n$ a concrete analytic $n$-gap in $\mathcal{P}(\omega)/fin$,  $\{\mathcal{I}_i : i\in n\}$ which is not weakly countably separated. Roughly speaking, the main result states that an analytic multiple gap $\{I_i : i\in n\}$ in $\mathcal{P}(\omega)/fin$ is either weakly countably separated or there is a copy of the object $\{\mathcal{I}_i : i\in n\}$ with $\mathcal{I}_i\subset I_i$. The aforementioned $n$-gap $\{\mathcal{I}_i : i\in n\}$ turns out to be a jigsaw. We also provide an example of an analytic multiple gap which is not a jigsaw. More results on the theory of analytic multiple gaps will appear in \cite{stronggaps} and further works.\\

In Section \ref{countablyseparated} we deal with multiple gaps which are weakly countably separated. In contrast with the results established in previous section, now given a weakly countably separated multiple gap $\{I_i : i \in n\}$ in $\mathcal{P}(\omega)/fin$ we find ideals of a particular form $\{I'_i : i\in n\}$ such that $I_i\subset I'_i$ (that is, above instead of below $I_i$). The most informative case is when the gap is indeed strongly countably separated, since then these new ideals above constitute a multiple gap. We show that every strongly countably separated multiple gap in $\mathcal{P}(\omega)/fin$ is below a strongly countably separated multiple gap which is moreover a dense jigsaw. In particular every dense strongly countably separated multiple gap in $\mathcal{P}(\omega)/fin$ is a jigsaw. We also provide an example of a strongly countably separated clover. This example is neither analytic nor dense.\\

In Section \ref{noHausdorff} we show that there is no analogue of Hausdorff's gap when we deal with $n$-gaps for $n\geq 3$. Namely, from Hausdorff's gap one gets that there exists a 2-gap $\{I_0,I_1\}$ in $\mathcal{P}(\omega)/fin$ such that $I_0$ and $I_1$ are generated by sets of size $\aleph_1$. We prove that if $MA_\theta$ holds for a cardinal $\theta$, then for $n\geq 3$ there exists no $n$-gap made of $\theta$-generated ideals in $\mathcal{P}(\omega)/fin$. Indeed we show that for every $n\geq 2$ it is consistent that there exists an $n$-gap of $\theta$-generated ideals in $\mathcal{P}(\omega)/fin$ but there are no such $n+1$-gaps. The technique used allows to prove some other results related to small sets both in $\mathcal{P}(\omega)/fin$ and in the Banach space $\ell_\infty/c_0$.\\

In Section \ref{completelyseparablesection} we construct, for every family $\frak{X}$ of subsets of $n$, a dense multiple gap which is a $B$-jigsaw for $B\in\frak{X}$ and a $B$-clover for $B\not\in\frak{X}$. The construction is based on the existence of completely separable almost dijoint families. It is an open problem whether such almost disjoint families exist in ZFC, though they are known to exist under several hypotheses.\\

In Section \ref{injectiveBanach}, we relate jigsaws with injectivity properties of Banach spaces. Let $\mathcal{B}$ be a Boolean algebra and $K$ its dual compact space in Stone's duality. We show that if $\mathcal{B}$ has a dense $n$-jigsaw for every $n$, then the Banach space of continuous functions $C(K)$ is uncomplemented in a superspace of the same density. Since we have found dense jigsaws in $\mathcal{P}(\omega)/fin$ of every size, this provides an alternative proof of a result of Amir \cite{Amir} that the Banach space $\ell_\infty/c_0$ is not injective.\\

We found that many notable examples of multiple gaps in $\mathcal{P}(\omega)/fin$ happen to be jigsaws, and jigsaws play an important role in the general theory. However we do not know that much about clovers. In particular we do not know whether analytic clovers -or any mixed version of clovers and jigsaws- exist at all, and we do not know if the results of Section~\ref{completelyseparablesection} hold in ZFC, except for the existence of dense jigsaws.\\

\section{Multiple gaps that are not countably separated}\label{dichotomy}

\subsection{An example of an analytic jigsaw that is not weakly countably separated}\label{analyticjigsaw}

Let $n<\omega$ and $T=n^{<\omega}$ be the $n$-ary tree. For every branch $x\in [T]$ and $i\in n$, we call $$a_x^i = \{s\in T : s^\frown i\in x\}.$$ Let $\mathcal{J}_i$ is the ideal of $\mathcal{P}(T)/fin$ generated by $\{a_x^i : x\in [T]\}$. Notice that each of these ideals is analytic, since its set of generators is an analytic set.

\begin{thm}\label{analyticjigsawexample}
The ideals $\{\mathcal{J}_i : i\in n\}$ constitute an $n$-jigsaw which is not weakly countably separated.
\end{thm}

Proof: It is obvious that the ideals are mutually orthogonal since $a_x^i\cap a_y^j =^\ast\emptyset$ whenever $i\neq j$. We prove first that they are not weakly countably separated, hence in particular they are a multiple gap. Let us denote by $T(m)$ the set of all elements of $t$ of length less than $m$. Suppose that the ideals were weakly countably separated, so that we have elements $c_i^k\subset T$, $i\in n$, $k<\omega$, such that $\bigcap_{i\in n}c_i^k = \emptyset$ for each $k$, and whenever we pick $b_i\in\mathcal{J}_i$ there exists $k$ with $b_i\subset^\ast c_i^k$ for every $i\in n$. In particular, for every $x\in [T]$ and every $i\in n$ there exist $k(x),m(x)\in \omega$ such that $a_x^i \setminus T(m(x)) \subset c_i^{k(x)}$. Let $X(m,k) = \{x\in [T] : \forall i\in n\ a_x^i \setminus T(m) \subset c_i^{k}\}$. Notice that each $X(m,k)$ is a closed subset of $[T]$ and that $\bigcup_{m,k<\omega}X(m,k) = [T]$. By Baire category theorem, there exist $m_0,k_0<\omega$ such that $X(m_0,k_0)$ has nonempty interior. This means that we can find $t\in T$ (that we can choose of length greater than $m_0$) such that $\{x\in [T] : t\in x\}\subset X(m_0,k_0)$. But this actually implies that $\{s\in T : s>t\}\subset\bigcap_{i\in n}c_i^{k_0}$, a contradiction.\\

We pass now to the proof that the ideals actually constitute a jigsaw.
Let $B\subset A\subset n$, and let $a\subset T$ such that $\{\mathcal{J}_i|_a : i\in A\}$ is a multiple gap. We have to find
 $b\subset a$, such that $b\perp \mathcal{J}_i$ for $i\in
n\setminus B$, while the ideals $\{\mathcal{J}_j|_b : j\in B\}$
are a multiple gap.\\

A nonempty subset $R\subset T$ is called a complete $A$-tree if for every $r\in R$ and every $i\in A$ there exists $s\in R$ such that $r^\frown i < s$.\\

\emph{Claim 1}: There exists a complete $A$-tree $R\subset a$. Proof of the claim: Let us call a set $u\subset T$ to be $i$-WF (for $i$-well founded), if it contains no infinite sequence $\{s_k : k<\omega\}$ with $s_k^\frown i < s_{k+1}$. Notice that an $i$-WF set is orthogonal to the ideal $\mathcal{J}_i$, therefore its complement is an upper bound of the ideal $\mathcal{J}_i$. It follows that $a$ cannot be a union of the form $a=\bigcup_{i\in A}u_i$ with each $u_i$ being $i$-WF. Otherwise, we could take $c_i = a\setminus u_i\geq \mathcal{J}_i|_a$ and $\bigcap c_i = \emptyset$, contradiciting that $\{\mathcal{J}_i|_a : i\in A\}$ is a multiple gap. We define a derivation procedure on subsets of $T$. Given a set $v\subset T$ we define $v[i]=\{s\in v : \not\exists t\in v : s^\frown i < t\}$ for $i\in A$, and then the derived set $v' = v\setminus \bigcup_{i\in A} v[i]$. By induction consider the iterated derivatives $a^{(\alpha)}$ of the set $a$ for ordinals $\alpha$. Choose an ordinal $\beta$ with $a^{(\beta)} = a^{(\beta+1)}$. This means that $a^{(\beta)}$ is a complete $A$-tree provided it is nonempty. But notice that $\bigcup_{\alpha<\beta}a^{(\alpha)}[i]$ is an $i$-WF set, so since $a$ cannot be the union of $i$-WF sets for $i\in A$, it follows that $a^{(\beta)}\neq\emptyset$. This finishes the proof of the claim.\\

Let $R$ be the complete $A$-tree given by Claim 1. We can find $b\subset R$ which is a complete $B$-tree and moreover $i$-WF for all $i\in n\setminus B$. We construct $b$ inductively, start with $b[0]\subset R$ a singleton, and then construct $b[n+1]$ by adding one element above $s^\frown j$ for every $s\in b[n]$ and $j\in B$. At the end, $b=\bigcup_{n<\omega}b[n]$. The set $b$ is orthogonal to $\mathcal{J}_i$, $i\not\in B$ since it is $i$-WF. So it remains to show that the ideals $\{\mathcal{J}_j|b : j\in B\}$ are a multiple gap. But this is the same proof as we already did, in the first paragraph of this proof, to show that the $\{\mathcal{J}_i : i\in n\}$ are a multiple gap, but for the integer $|B|$ instead of $n$. The reason is that we can identify the set $b$ with the $|B|$-adic tree $B^{<\omega}$ and then the ideals $\mathcal{J}_i|b$ are identified with the originally considered ideals $\mathcal{J}_i$.$\qed$\\

\subsection{The dichotomy}

We show that if an analytic multiple gap is not weakly countably separated, then it contains a copy the Example \ref{analyticjigsaw} $\{\mathcal{J}_i : i\in n\}$ of analytic jigsaw presented above.

\begin{thm}[Analytic multiple gap dychotomy]
Let $\{I_i : i\in n\}$ be analytic ideals of $\mathcal{P}(\omega)$ which constitute a multiple gap of $\mathcal{P}(\omega)/fin$. Then,
\begin{enumerate}
\item either the ideals are weakly countably separated in $\mathcal{P}(\omega)/fin$,
\item or there exists a one-to-one function $u:T\To \omega$, where $T=n^{<\omega}$ is the $n$-ary tree, such that $u(\mathcal{J}_i)\subset I_i$ for all $i\in n$.
\end{enumerate}

\end{thm}

Proof: Let $\Upsilon = \omega^{<\omega}$ be the Baire tree. Let $\Sigma = \{(e_i)_{i\in n}\in \prod_{i\in n}I_i : e_i\cap e_j=\emptyset\text{ for }i\neq j\}$. Since $\Sigma$ is an analytic set, there exists a continuous surjection $f:[\Upsilon]\To \Sigma$. We write $f(x) = (f_i(x))_{i\in n}$. We distinguish two cases.\\

In the first case, we suppose that there exists a countable decomposition $[\Upsilon] = \bigcup_{k<\omega}X_k$ such that

$(\star)$ for every $k$ and every $\{x_i : i\in n\}\subset X_k$ we have that $\bigcap_{i\in n}f_i(x_i)=\emptyset$.

In this case we are in alternative (1) of the dichotomy. The elements $c_i^k = \bigcup_{x\in X_k} f_i(x)$ witness countable separation.\\

For the second case, we assume that the first case does not hold. Let $\Upsilon'\subset\Upsilon$ be the family of all $s\in\Upsilon$ such that $[\Upsilon_s]$ cannot be decomposed into countably many pieces $[\Upsilon_s] = \bigcup_{k<\omega}X_k$ with the property $(\star)$ specified in the first case. Notice that $\Upsilon'$ is closed under initial segments.\\

\emph{Claim 1.} If $t\in\Upsilon'$, then $[\Upsilon'_t]$ cannot be decomposed into countably many pieces with property $(\star)$. Proof of the claim: Consider $F$ the set of all $s>t$ such that $s\not\in\Upsilon'$. Then $[\Upsilon_t] = [\Upsilon'_t]\cup \bigcup_{s\in F}[\Upsilon_s]$. All the sets $[\Upsilon_s]$ in the right admit countable decompositions as $(\star)$, so since $[\Upsilon_t]$ does not have such a decomposition, neither has $[\Upsilon'_t]$.\\

\emph{Claim 2.} If $t\in\Upsilon'$, and $F\subset\omega$ is a finite set, then there exist branches $\{x_i : i\in n\}\subset [\Upsilon'_t]$ such that $\bigcap_{i\in n}f_i(x_i)\setminus F\neq\emptyset$. Proof of claim 2. Suppose that there existed a finite set $F$ such that $\bigcap_{i\in n}f_i(x_i)\subset F$ for all $\{x_i : i\in n\}\subset [\Upsilon'_t]$. Then, we can find a countable (even finite) decomposition $[\Upsilon'_t] = \bigcup_{m<\omega} X_m$ and finite sets $F_i(m)\subset F$ for $i\in n$ and $m\in\omega$ such that $f_i(x)\cap F = F_i(m)$ for every $x\in X_m$. Since $\bigcap_{i\in n}f_i(x) = \emptyset$ for every $x$, it follows that $\bigcap_{i\in n} F_i(m) = \emptyset$ for every $m<\omega$. Hence for every $\{x_i : i\in n\}\subset X_m$, $\bigcap_{i\in n}f_i(x_i)  = \emptyset$. This means that every $X_m$ has property $(\star)$, in contradiction with Claim 1.\\

\emph{Claim 3.} Let $s\in\Upsilon'$ and $F\subset \omega$ finite. Then there exist $\{s_i : i\in n\}\subset \Upsilon'$ incomparable elements above $s$ and $k\in\omega\setminus F$ such that $k\in \bigcap_{i\in n}f_i(x_i)$ whenever $s_i\in x_i$ for every $i$. Proof of the claim: By Claim 2, there exist $\{y_i : i\in n\}\subset [\Upsilon'_s]$ such that $\bigcap_{i\in n}f_i(y_i)\setminus F\neq\emptyset$. Fix $k\in \bigcap_{i\in n}f_i(y_i)\setminus F$. Now, using the fact that $f$ is continuous we can find elements $s_i\in\Upsilon'_t$ such that $s<s_i\in y_i$ and which satisfy the requirement of the claim.\\

We construct inductively two functions $v:T\To \Upsilon'$ and $u:T\To \omega$. $v(\emptyset)$ is just the root of $\Upsilon'$. Once $v(t)$ is defined, we apply Claim 3 to $s=v(t)$ and $F$ the set of all $u(t')$'s that have been previously chosen in the inductive procedure (in order to make $u$ one-to-one). We set $v(t^\frown i) = s_i$ and $u(t) = k$ given by Claim 3. In this way, we have the property that for every $t\in T$, if $x\in [\Upsilon']$ and $v(t^\frown i)\in x$, then $u(t)\in f_i(x)$. To conclude the proof, we have to check that $u(\mathcal{J}_i) \subset I_i$. So let $x\in [T]$ and $a=a_i^x = \{t\in T : t^\frown i \in x\}$ one of the generators of $\mathcal{J}_i$. We want to see that $u(a)\in I_i$. The set $\{v(t) : t\in x\}$ is an infinite chain in $\Upsilon'$, hence cofinal in some branch $y\in [\Upsilon']$. If $t\in a$, then $v(t^\frown i)\in y$, hence $u(t)\in f_i(y)$. It follows that $u(a)\subset f_i(y) \in I_i$. $\qed$\\

\subsection{Example of an analytic gap that is not a jigsaw} Let $T=2^{<\omega}$ be the dyadic tree. Again, for $i=0,1$, we define elements $a_x^i = \{s\in T : s^\frown i\in x\}$ and the ideal $I_i$ of $\mathcal{P}(T)/fin$ as generated by all $a_x^i$, for $x\in [T]$. On the other hand, let $J$ be the ideal generated by the antichains of $T$.

\begin{prop}
The ideals $\{I_0,I_1,J\}$ are three analytic ideals which are a multiple gap which is not weakly countably separated and that is neither a jigsaw nor a clover.
\end{prop}

Proof: The fact that they are analytic is easily checked as they have analytic sets of generators. We see that they are not weakly countably separated (in particular, they are a 3-gap). Let us call $J=I_2$ for convenience. Suppose we have $c_i^k$ with $\bigcap_{i\in 3}c_i^k =^\ast\emptyset$ for every $k<\omega$, and for every $x_i\in  I_i$ ($i\in 3$) there is $k<\omega$ with $x_i\leq c_i^k$. The proof is similar to that of Theorem~\ref{analyticjigsawexample}. We call $T(m)$ to the set of elements of $T$ of length less than $m$. For every $x\in [T]$ there exist $k(x)$ and $m(x)$ integers such that $a_x^i \subset c_i^{k(x)} \setminus T(m(x))$ for $i=0,1$ and $$\{s^\frown 0 : s\in a_x^1\}\cup\{s^\frown 1 : s\in a_x^0\}\subset c_2^{k(x)}\setminus T(m(x)).$$
By a Baire category argument, there exist $k,m<\omega$ and $t\in T$ of length greater than $m$ such that the properties above hold for all $x\in [T_t]$ taking $k(x)=k$ and $m(x)=m$. It follows that $t'\in\bigcap_{i\in 3}c_i^k$ for every $t'>t$. This shows that $\{I_0,I_1,J\}$ are not weakly countably separated.

\emph{Claim 1}: if $b\in J^\perp$, then $b$ is contained in a finite number of chains. Proof of the claim: Let $A = \{x\in [T] : b\cap x \neq^\ast \emptyset\}$. Notice that $b\setminus \bigcup A$ is an antichain, so $b\subset^\ast \bigcup A$. If we suppose that $A$ is infinite, then it contains either an increasing or a decreasing sequence $\{x_n\}$ in the lexicographical order of $[T]=2^\omega$. If $t_n\in x_n\cap b$ is chosen high enough so that $t_n> x_n\wedge x_{n+1}$, then $\{t_n : n<\omega\}\subset b$ is an antichain. This finishes the proof of the claim. 

The 3-gap is not a jigsaw, because if $b\in J^\perp$, then by claim 1 above, $I_0|_b$ and $I_1|_b$ are finitely generated and therefore separated. The 3-gap is neither a clover. Consider $c$ the set of all $(k_0,k_1,\ldots,k_{m-1})\in T$ such that $m$ is even and $k_j=1$ for all $k$ even. Then $c$ is a dyadic subtree of $T$, $c\in I_0^\perp$ and $\{I_1|_c,J|_c\}$ are a gap (this follows from Claim 1, because $I_1|_c$ is the ideal generated by the chains of $c$ and $I_2|_c$ the ideal generated by the antichains).$\qed$\\

\section{Multiple gaps that are countably separated}\label{countablyseparated}

Let $D$ be a countable subset of a compact metric space $L$ and let $K = acc(D)$ be the set of accumulation points of $D$, those points such that every neighborhood meets infinitely many elements of $D$. For a subset $G\subset K$, we define

$$\mathcal{I}_G = \{a\subset D : acc(a)\subset G\}.$$

Note that $\mathcal{I}_G$ is an ideal of $\mathcal{P}(D)/fin$. Note also the fact that $\mathcal{I}_{G}$ is orthogonal to $\mathcal{I}_{G'}$ if and only if $G\cap G' = \emptyset$. Theorem~\ref{caracterizationcs} below shows that this kind of ideals characterize countable separation. 

\begin{lem}\label{Xseparation}
Let $L$ be a compact space, $\mathfrak{X}$ a family of subsets of $n$, $\mathfrak B$ a basis for the topology of $L$ which is closed under finite unions and intersections, and $\{Z_i : i\in n\}$ closed subsets of $L$ such that $\bigcap_{i\in A}Z_i = \emptyset$ whenever $A\in \mathfrak X$. Then there exist sets $Z_i\subset V_i\in\mathfrak B$ with $\bigcap_{i\in A}\overline{V_i} = \emptyset$ whenever $A\in\mathfrak X$.
\end{lem}

Proof:  First we prove the case in which $\mathfrak X$ consists only of the full subset of $n$. This is proved by induction on $m$, the case $m=2$ being just the standard compactness argument for separation. For $m>2$, $Z = \bigcap_{i\in m-1}Z_i$ and $Z_{m-1}$ are disjoint closed subsets of $L$, so by case $m=2$, there are $V,W\in\mathfrak B$ with $Z_{m-1}\subset V$, $Z\subset W$ and $\overline{V}\cap \overline{W} = \emptyset$. The inductive hypothesis provides $\{U_i : i\in m-1\}\subset\mathfrak B$ with $\bigcap_{i\in m-1}\overline{U_i} = \emptyset$ and $Z_i\setminus W\subset U_i$ for $i\in m-1$, and then we define $V_{m-1} = V$ and $V_i = U_i\cup W$ for $i\in m-1$. Now, for a general family $\mathfrak X$ of subsets of $n$, the previous case allows to find, for every $A$, $Z_i\subset V_i^A\in \mathfrak B$ such that $\bigcap_{i\in A}\overline{V_i^A} = \emptyset$, and take $V_i = \bigcap_{A\in\mathfrak{X}}V_i^A$.$\qed$\\

\begin{thm}\label{caracterizationcs}
Let $\{I_i : i\in n\}$ be ideals of $\mathcal{P}(\omega)/fin$ and let $\mathfrak{X}$ a family of subsets of $n$. The following are equivalent:

\begin{enumerate}
\item The ideals are $\mathfrak{X}$-countably separated.
\item There exists a continuous function\footnote{We recall that we denote by $\omega^\ast$ the Stone space of $\mathcal{P}(\omega)/fin$} $\phi:\omega^\ast\To 2^\omega$ such that $\bigcap_{i\in A}\phi(U(I_i)) = \emptyset$ for every $A\in\mathfrak{X}$.
\item There exists a bijection $f:\omega\To D$, where $D$ is a countable subset of a compact metric space $L$, and sets $G_i\subset K = acc(D)$, $i\in n$ such that
\begin{enumerate}
\item $\bigcap_{i\in A}G_i = \emptyset$ for every $A\in\mathfrak{X}$, and
\item $f(I_i) = \{f(a) : a\in I_i\}\subset \mathcal{I}_{G_i}$ for every $i\in n$.
\end{enumerate}
\item The previous condition holds for $D = \{\xi\in 2^\omega : \exists m_0 \forall m>m_0\ \xi_m=0\}$.
\end{enumerate}

\end{thm}

Proof: We prove $4\Rightarrow 3\Rightarrow 1\Rightarrow 2\Rightarrow 4$. First implication is obvious, so we start by showing that $3\Rightarrow 1$. We prove that the ideals $\{\mathcal{I}_{G_i} : i\in n\}$ are $\mathfrak{X}$-countably separated. Let $\mathfrak{B}$ be a countable basis for the topology of $L$ which is closed under finite unions and intersections. Consider the family of all $n$-tuples of subsets of $D$ of the form $(V_i\cap D)_{i\in n}$ such that $V_i\in\mathfrak{B}$ and $\bigcap_{i\in A}\overline{V_i}=\emptyset$ for all $A\in\mathfrak{X}$. This is a countable family of $n$-tuples of elements of $\mathcal{P}(D)$. We see that it witnesses $\mathfrak{X}$-countable separation. The first condition to check is obvious: $\bigcap_{i\in A}V_i\cap D = ^\ast \emptyset$ for $A\in\mathfrak{X}$, whenever $(V_i\cap D)_{i\in n}$ is in the family. For the second condition, we pick $a_i\in\mathcal{I}_{G_i}$. Since $acc(a_i)\subset G_i$ for each $i$, we have that $\bigcap_{i\in A}acc(a_i) = \emptyset$ for $A\in\mathfrak{X}$. By Lemma~\ref{Xseparation}, we can find $V_i\in\mathfrak{B}$ such that $acc(a_i)\subset V_i$ and $\bigcap_{i\in A}\overline{V_i} = \emptyset$ for all $A\in\mathfrak{X}$. It follows that $a_i\subset^\ast V_i\cap D$ and the tuple $(V_i\cap D)_{i\in n}$ belongs to our family.

We pass now to implication $1\Rightarrow 2$. Using Stone's duality, $\mathfrak{X}$-countable separation can be restated as the existence of clopen sets $\{C_m^i: i\in n, m\in\omega\}$ of $\omega^\ast$ such that
\begin{itemize} \item $\cap_{i\in A}C_m^i = \emptyset$ for every $m\in\omega$ and every $A\in\mathfrak{X}$, and
\item for every clopen sets $x_i\subset U(I_i)$. $i\in n$, there exists $m$ such that $x_i\subset C_m^i$ for every $i\in n$.
\end{itemize}

For every $(m,i)\in\omega\times n$ define $\phi_{(m,i)}:\omega^\ast\To 2$ as the characteristic function of the clopen set $C_m^i$. All these functions together provide a continuous function $\phi:\omega^\ast\To 2^{\omega\times n}$. Suppose now for contradiction that we had $\xi\in\bigcap_{i\in A}\phi(U(I_i))$ with $A\in\mathfrak{X}$. Then for every $i\in A$, $\xi = \phi(t_i)$ with $t_i\in U(I_i)$. There must exist $m\in\omega$ such that $t_i\in C_m^i$ for every $i\in A$. For every $i,j\in A$ we have that
$$1 = \phi_{(m,i)}(t_i) = \xi_{(m,i)} = \phi_{(m,i)}(t_j).$$
It follows that $t_j\in\bigcap_{i\in A}C_m^i = \emptyset$, a contradiction.

We prove now that $2\Rightarrow 4$. Given $\phi$, we construct inductively a tree of subsets of $\omega$, $\{a_s : s\in 2^{<\omega}\}$ with
\begin{enumerate}
\item $\phi^{-1}(\xi\in 2^\omega : \xi|_{length(s)} = s) = clopen(a_s)$, where $clopen(a)$ denotes the clopen subset of $\omega^\ast$ associated to the set $a\in\mathcal{P}(\omega)$.
\item $a_\emptyset = \omega$,
\item $a_s = a_{s^\frown 0}\cup a_{s^\frown 1}$,
\item $\emptyset = a_{s^\frown 0}\cap a_{s^\frown 1}$,
\item If we call $m_s = \min(a_s\setminus\{m_t : t< s\})$, then $m_s\in a_r$ whenever $s<r$ and $r_k=0$ for all $k\geq length(s)$.
\end{enumerate}

Notice that the function $s\mapsto m_s$ is a bijection from $2^{<\omega}$ onto $\omega$. Consider $g:\omega\To 2^{<\omega}$ its inverse function, so that $m_{g(m)} = m$ for all $m$. We define $f:\omega\To D$ as $f(m) = g(m)^\frown 1^\frown(0,0,0,\ldots)$. Let $G_i = \phi(U(I_i))$. We have to prove that if $x\in I_i$, then $acc\{f(m) : m\in x\}\subset G_i$. So fix $x\in I_i$ and suppose that $\xi\in acc\{f(m) : m\in x\}$. There exists $\{m_1,m_2,\ldots\}\subset x$ such that the sequence $\{f(m_k)\}_{k<\omega}$ converges to $\xi$ in the space $2^\omega$. This means that for every $k\in\omega$ all but finitely many $p\in\omega$ satisfy $(\xi_0,\ldots,\xi_k)< g(m_p)$. This implies that, for every $k$,

$$\{m_0,m_1,\ldots\}\subset^\ast a_{(\xi_0,\ldots,\xi_k)}.$$

This is because if $t<g(m)$ then $m\in a_t$ (notice that $m_s\in a_s$, which implies $m = m_{g(m)}\in a_{g(m)}\subset a_t$). Now, by Stone's duality

\begin{eqnarray*}
clopen\left\{m_0,m_1,\ldots\right\} &\subset & clopen\left(a_{(\xi_0,\ldots,\xi_k)}\right)\\
&=& \phi^{-1}\{\zeta\in 2^\omega : \zeta_i = \xi_i, i=0,\ldots,k\}
\end{eqnarray*}

It follows that $clopen\{m_0,m_1,\ldots\}\subset \phi^{-1}(\xi)$. But $$clopen\{m_0,m_1,\ldots\} \subset clopen(x),$$ and $clopen(x)\subset U(I_i)$ since $x\in I_i$. Hence $\xi\in\phi(U(I_i)) = G_i$ as we wanted to prove.$\qed$\\

In Theorem \ref{caracterizationcs}, notice that if the ideals $\{I_i : i\in n\}$ are strongly countably separated, then the sets $G_i$ are pairwise disjoint, and then the ideals $\{\mathcal{I}_{G_i} : i\in n\}$ are mutually orthogonal. Hence, in that case, if $\{I_i : i\in n\}$ are a multiple gap, then the larger ideals $\{\mathcal{I}_{G_i} : i\in n\}$ are also a multiple gap.\\

\begin{prop}
In Theorem \ref{caracterizationcs}, if the ideals $I_i$ are analytic, then the sets $G_i$ can be taken to be Borel.
\end{prop}

Proof: Let us compute the complexity of the sets $G_i$ obtained in the proof of $2\Rightarrow 4$ of Theorem~\ref{caracterizationcs}. The third equality below follows from the fact that nonempty $G_\delta$ subsets of $\omega^\ast$ have nonempty interior.

\begin{eqnarray*}
G_i &=& \phi(U(I_i)) = \{\xi\in 2^\omega : \exists t\in U(I_i)\ \phi(t) = \xi\}\\
&=& \{\xi\in 2^\omega : \exists c\subset U(I_i),\text{ c clopen, }\phi|_c = \xi\}\\
&=& \{\xi\in 2^\omega : \exists a\in I_i\ \phi|_{clopen(a)} = \xi\}\\
&=& \{\xi\in 2^\omega : \exists a\in I_i\ \forall k\ \exists m,\ a\setminus m\subset a_{(\xi_0,\ldots,\xi_k)}\}.
\end{eqnarray*}

Hence, if the ideals $I_i$ are analytic, then the sets $G_i$ are analytic. Using a generalized version of Lusin's separation theorem (which is proved in the same manner as Lemma~\ref{Xseparation}), we can find Borel sets $G'_i$ with $G_i\subset G'_i$ and still $\bigcap_{i\in A}G'_i = \emptyset$ for all $A\in\mathfrak{X}$.$\qed$\\

Notice that the fact that $G$ is Borel does not mean that the ideal $\mathcal{I}_{G}$ is Borel. Indeed one can check that $\mathcal{I}_{\mathbb{Q}}$ is true coanalytic. However, when $G$ is $G_\delta$ we have:\\

\begin{prop}\label{IGBorel}
If $G$ is a $G_\delta$ set, then $\mathcal{I}_G$ is Borel, indeed $F_{\sigma\delta}$ set.
\end{prop}

Proof: Let $G=\bigcap_{m<\omega}U_m$ with $U_m$ open. For every $m$, let $U_m = \bigcup_{k<\omega}V^k_m$ where each $V^k_m$ is open, $V^1_m\subset V^2_m\cdots$ and $\overline{V^k_m}\subset U_m$. Then,
$$\mathcal{I}_G = \bigcap_{m<\omega}\bigcup_{k<\omega}\bigcup_{F\subset\omega\ finite}\{a\subset D : a\setminus F\subset V_m^k\}$$
and the sets that appear in the right are closed. Let us check each inclusion. For $[\subset]$, assume that $a\in \mathcal{I}_G$. Then $acc(a)\subset U_m$ for every $m$. Fix such an $m$. By compactness of $acc(a)$, there exists $k$ such that $acc(a)\subset V^k_m$. Since $V^k_m$ is open and $L$ is compact, it follows that there exists $F\subset \omega$ finite such that $a\setminus F\subset V^k_m$. For $[\supset]$, assume that
for every $m$ there exist $k$ and a finite $F\subset \omega$ such that $a\setminus F\subset V_m^k$. Then $acc(a)\subset \overline{V_m^k}\subset U_m$ for every $m$. If follows that $acc(a)\subset G$.
$\qed$\\

For the rest of this section, $L$ is a compact metric space, $D\subset L$ is countable and $K=acc(D)$.\\

\begin{lem}\label{relativeopen}
Let $Z\subset L$ and let $\{V_i : i\in n\}$ be relative open subsets of $Z$ with $\bigcap_{i\in n}V_i = \emptyset$. Then there exist $\{U_i : i\in n\}$ open subsets of $L$ such that $U_i\cap Z = V_i$ and $\bigcap_{i\in n}U_i=\emptyset$.\\
\end{lem}

Proof: It is enough to prove the case when $Z$ is closed in $L$ (if $Z$ is not closed, we can deal with relative open sets of $\overline{Z}$). Let $f_i:Z\To [0,1]$ be continuous functions such that $f_i^{-1}(0) = Z\setminus V_i$. Let $f:Z\To [0,1]^n$ be given by putting together the functions $f_i$, $i<n$. Then $f(Z)\subset H = \{(t_i)_{i\in n} : \exists i\ t_i = 0\}\subset [0,1]^n$. Notice that $H$ is a retract of $[0,1]^n$, the retraction being $r((t_i)_{i\in n}) = (t_i - \min[t_j : j\in n])_{i\in n}$. Hence, there exists a continuous function $\hat{f}:L\To H$ such that $\hat{f}|_Z = f$. The sets $U_j = \hat{f}^{-1}((t_i)_{i\in n}\in H : t_j>0)$ are the ones we are looking for.$\qed$\\

Given a countable compact set $S$, we denote by $S'$ its Cantor-Bendixson derivative, that is, the set of all points of $S$ which are not isolated in $S$. Inductively, we define $S^{(n)} = [S^{(n-1)}]'$, assuming that $S^{(0)} = S$. The space $S$ has height $n$ if $S^{(n)}=\emptyset$ but $S^{(n-1)}\neq \emptyset$.

\begin{thm}\label{Giaregap}
Let $\{G_i : i\in n\}$ be pairwise disjoint subsets of $K$. The following are equivalent
\begin{enumerate}
\item $\{\mathcal{I}_{G_i} : i\in n\}$ is an $n$-gap.
\item For every open subsets $\{U_i : i\in n\}$ of $L$ such that $G_i\subset U_i$ for every $i\in n$, we have that $\bigcap_{i\in n}U_i \neq \emptyset$.
\item There exists a countable compact $S\subset K$ of height $n$, and a bijection $\sigma:n\To n$ such that
\begin{enumerate}
\item $S^{(j)}\setminus S^{(j+1)}\subset G_{\sigma(j)}$ for $j\in n$.
\end{enumerate}
\end{enumerate}
\end{thm}

Proof: For $1\Rightarrow 2$, suppose that there exist open sets $U_i\subset L$, $G_i\subset U_i$, $\bigcap_{i\in n}U_i = \emptyset$. Let $a_i = D\cap U_i$. Then $\bigcap_{i\in n}a_i = \emptyset$ and if $x\in \mathcal{I}_{G_i}$, then $x\subset^\ast a_i$. Hence the ideals $\{\mathcal{I}_{G_i} : i\in n\}$ are separated by the $a_i$'s and they are not an $n$-gap. For the converse implication $2\Rightarrow 1$, suppose that the ideals are not a multiple gap, so there exist $a_i\subset D$ such that $\bigcap_{i\in n}a_i = \emptyset$ and $x\subset^\ast \mathcal{I}_{G_i}$ whenever $x\in \mathcal{I}_{G_i}$. Let $V_i = K\setminus acc(D\setminus a_i)$. These are relative open subsets of $K$. On the one hand,
\begin{eqnarray*}
\bigcap_{i\in n}V_i &=& K\setminus \bigcup_{i\in n}acc(D\setminus a_i) = K\setminus acc\left(\bigcup_{i\in n}D\setminus a_i\right)\\ &=& K\setminus acc\left(D\setminus \bigcap_{i\in n} a_i\right) = K\setminus acc(D) = \emptyset
\end{eqnarray*}

On the other hand, $G_i\subset V_i$ because if $\xi\not\in V_i$, then there is a convergent sequence $(d_n)$ contained in $D\setminus a_i$ that converges to $\xi$. If, in addition $\xi\in G_i$, then $x = \{d_n : n<\omega\}\in \mathcal{I}_{G_i}$, hence $x\subset^\ast a_i$, a contradiction. We use Lemma~\ref{relativeopen} to get the open sets $U_i$ from the relative open sets $V_i$.\\

For $3\Rightarrow 2$, we proceed by induction on $n$. Suppose that we have $G_i\subset U_i$ with $U_i$ open, $i\in n$. Then $G_{\sigma(n-1)}\subset U_{\sigma(n-1)}$, hence $U_{\sigma(n-1)}$ contains a point $\xi_1$ of the level $S^{(n-1)}$ of the scattered space $S$. Being open, this implies that $U_{\sigma(n-1)}$ contains a copy of a countable compact space $R$ of height $n-1$ contained in $S\setminus S^{(n-1)}$ which is a clopen of $S$. Notice that $R^{(j)}\setminus R^{(j+1)} = (S^{(j)}\setminus S^{(j+1)})\cap R$ for every $j\in n-1$. Hence $R^{(k)}\setminus R^{(k+1)} \subset U_{\sigma(k)}\cap U_{\sigma(n-1)}$ for $k\in n-1$. We can apply then the inductive hypothesis to conclude that $\bigcap_{k\in n-1} U_{\sigma(k)}\cap U_{\sigma(n-1)}\neq\emptyset$.\\

The implication that requires more work is $2\Rightarrow 3$. We suppose that $(3)$ does not hold and we will find open sets $V_i\supset G_i$ with $\bigcap_{i\in n}V_i = \emptyset$. We proceed by induction on $n$. In case $n=2$, the failure of condition $(3)$ simply means that there is no sequence of elements of $G_0$ converging to a point of $G_1$ nor viceversa. Then $G_0\subset V_0 = \overline{G_0\cup G_1}\setminus\overline{G_1}$ and $G_1\subset V_1 = \overline{G_0\cup G_1}\setminus\overline{G_0}$. The sets $V_0$ and $V_1$ are relative open in $\overline{G_0\cup G_1}$ and $V_0\cap V_1=\emptyset$, hence by Lemma~\ref{relativeopen}, $G_0$ and $G_1$ are separated by open sets. Now we prove the implication for a fixed integer $n$ assuming it holds for all integers less than $n$. We start with a claim showing that separation by open sets holds locally:\\

Claim 1: If $\xi\in \bigcup_{i\in n}G_i$ then there exists an open neighborhood $W$ of $\xi$ and open sets $V_j\subset W$ for $j\in n$, such that  $G_j\cap W\subset V_j$ and $\bigcap_{j\neq n}V_j = \emptyset$. Proof of the claim: Fix $i\in n$ such that $\xi\in G_i$. We will prove that the claim holds with the additional fact that $V_i= W$. If, reasoning towards contradiction, we were unable to find $\{V_j : j\neq i\}$ as required for any neighborhood $W$ of $\xi$, then the sets $\{G_j\cap W : j\neq i\}$ cannot be separated by open sets for any neighborhood $W$ of $\xi$. By inductive hypothesis, this means that we can find a bijection $\sigma:n-1\To n\setminus\{i\}$ and a countable compact $S$ of height $n-1$ such that $S^{(k)}\setminus S^{(k+1)}\subset G_{\sigma(k)}\cap W$ for every $k\in n-1$. Since we can do this for every neighborhood of $\xi$, we can find a sequence $\{S_m : m<\omega\}$ of pairwise disjoint\footnote{Inductively, we find $S_m$ inside $W_m\setminus \bigcup_{m'<m}S_m$,where $\{W_m : m<\omega\}$ is a neighborhood basis of $\xi$} countable compact spaces of height $n-1$ such that
\begin{enumerate}
\item $S_m^{(k)}\setminus S_m^{(k-1)}\subset G_{\sigma(k)}$ for every $k\in n-1$ and every $m<\omega$;
\item $\overline{\bigcup_{m<\omega}S_m} = \{\xi\} \cup \bigcup_{m<\omega} S_m$.
\end{enumerate}
Now, the space $S_\infty = \{\xi\} \cup\bigcup_{m<\omega} S_m$ is a countable compact space such that $S_\infty^{(k)}\setminus S_\infty^{(k+1)}\subset G_{\sigma(k)}$ and $S_\infty^{(n-1)}\setminus S_\infty^{(n)} = \{\xi\} \subset G_{i}$. This means that condition (3) holds, and we are supposing that it does not. This finishes the proof of the claim.\\

Now let $\mathcal{F}_0$ be the family of all open subsets $W$ of $K$ such that there exists open sets $V_i(W)\subset W$ with $\bigcap_{i\in n}V_i(W) = \emptyset$ and $G_i\cap W\subset V_i(W)$. Let $\mathcal{F}_1$ be the family of all open subsets $W_1$ of $K$ such that there exists $W_0\in\mathcal{F}_0$ such that $\overline{W_1}\subset W_0$. Let $Z = \bigcup\mathcal{F}_0= \bigcup\mathcal{F}_1$. By Claim 1 above, $\bigcup_{i\in n} G_i\subset Z$. Since $Z$ is a metric space, it is paracompact~\cite{Stone}, so there exists $\mathcal{F}_2$ a locally finite open refinement of $\mathcal{F}_1$ which covers $Z$. For every $W\in \mathcal{F}_2$, we fix $W_0\in\mathcal{F}_0$ such that $\overline{W}\subset W_0$. For every $i\in n$, let $$V_i = \bigcap \{V_i(W_0)\cup (Z\setminus\overline{W}) : W\in\mathcal{F}_2\}.$$
We have the following properties:
\begin{enumerate}
\item Each $V_i$ is an open set because it is locally a finite intersection of open sets. Namely, since $\mathcal{F}_2$ is locally finite every point $\xi$ has a neighborhood $U\subset Z$ that meets only finitely many sets $W\in \mathcal{F}_2$. Let $U_0$ be an open neighborhood of $\xi$ with $\overline{U_0}\subset U$. Then $U_0\subset Z\setminus \overline{W}$ whenever $W\cap U = \emptyset$, hence
$$U_0\cap V_i = U_0\cap \bigcap\{V_i(W_0)\cup (Z\setminus\overline{W}) : W\cap U\neq\emptyset\}$$ is open because it is a finite intersection of open sets.
\item $G_i\subset V_i$ for every $i\in n$. Namely, if $\xi\in G_i$ and $W\in \mathcal{F}_2$ then we have two options. First case is that $\xi\in W_0$, and then $\xi\in G_i\cap W_0 \subset V_i(W_0)$. Second case is that $\xi\not\in W_0$, and then $\xi\in Z\setminus\overline{W}$ since $\overline{W}\subset W_0$.
\item $\bigcap_{i\in n}V_i = \emptyset$. Namely, suppose $\xi\in \bigcap_{i\in n}V_i$. Since $\xi\in Z$, consider $W\in\mathcal{F}_2$ such that $\xi\in W$. Then $\xi\not\in Z\setminus\overline{W}$, so since $\xi\in \bigcap_{i\in n}V_i$, it follows that $\xi\in\bigcap_{i\in n}V_i(W_0) = \emptyset$, a contradiction.
\end{enumerate}
This together with Lemma~\ref{relativeopen} finishes the proof.$\qed$\\

\begin{thm}\label{Giarejigsaw}
Let $\{G_i : i\in n\}$ be disjoint subsets of $K$ such that $\{\mathcal{I}_{G_i} : i\in n\}$ are an $n$-gap. Then the ideals $\{\mathcal{I}_{G_i} : i\in n\}$ are a jigsaw. Moreover, in the definition of jigsaw, the element $b$ can always be chosen so that the ideals $\mathcal{I}_{G_i}|_b$ are Borel $F_{\sigma\delta}$.
\end{thm}

Proof: Let $a\subset D$  and $A\subset n$ such that $\{\mathcal{I}_{G_i}|_a : i\in A\}$ form a multiple gap. Let $B\subset A$ and we have to find $b\subset a$ such that $\{\mathcal{I}_{G_i}|_b : i\in B\}$ is a multiple gap and $b\in\mathcal{I}_{G_i}^\perp$ for $i\in A\setminus B$. We can suppose, without loss of generality, that $D=a$, $acc(a)=K$ and $A=n$, because $\{\mathcal{I}_{G_i}|_a : i\in A\}$ is just the same as $\{\mathcal{I}_{G_i\cap acc(a)} : i\in A\}$ (taking in the latter case $a$ instead of $D$ as the ground countable set). By Theorem~\ref{Giaregap} there exists a countable compact $S\subset K$ of height $n$ and a bijection $\sigma:n\To n$ such that $S^{(j)}\setminus S^{(j+1)}\subset G_{\sigma(j)}$ for all $j\in n$. Let us enumerate the set $B$ as $B = \{\sigma(i_j) : j<r\}$ with $i_0<i_1<\cdots$. Elementary manipulation of countable compact spaces shows that we can find a compact subset $R\subset S$ of height $r$ such that $$R^{(j)}\setminus R^{(j+1)}\subset S^{(i_j)}\setminus S^{(i_j+1)}\subset G_{\sigma(i_j)}$$
for all $j\in r$. Also, it is possible to find $b\subset D$ such that $acc(b) = R$ (for example consider a distance $\rho$ metrizing $L$, enumerate $R=\{\xi_p : p<\omega\}$, let $\xi_p^m\in D$ such that $\rho(\xi_p^m,\xi_p)<\frac{1}{m}$ and set $b=\{\xi_p^m : m>p\}$). The set $b$ is the one we are looking for. On the one hand, $$acc(b) \subset R \subset \bigcup_{j<r}G_{\sigma(i_j)} = \bigcup_{i\in B}G_i,$$ hence $b\in\mathcal{I}_{G_i}^\perp$ for $i\not\in B$. On the other hand $\{\mathcal{I}_{G_i}|_b : i\in B\}$ is just the same as $\{\mathcal{I}_{G_i\cap R} : i\in B\}$ and this is a multiple gap by Theorem~\ref{Giaregap}, since $R^{(j)}\setminus R^{(j+1)}\subset G_{\sigma(i_j)}$. Moreover, the ideals $\mathcal{I}_{G_i\cap R}$ are Borel $F_{\sigma\delta}$ sets by Proposition~\ref{IGBorel}: since $R$ is countable, all its subsets are $G_\delta$.$\qed$\\

\begin{cor}\label{containedindensejigsaw}
If $\{I_i : i\in n\}$ is a strongly countably separated multiple gap in $\mathcal{P}(\omega)/fin$, then there exists a strongly countably separated dense jigsaw $\{J_i : i\in n\}$ such that $I_i\subset J_i$ for every $i\in n$.
\end{cor}

Proof: By Theorem~\ref{caracterizationcs} we have $I_i\subset\mathcal{I}_{G_i}$ where the $G_i$'s are pairwise disjoint. We can suppose that $\bigcup G_i = K$, so that $\{\mathcal{I}_{G_i} : i\in n\}$ is dense. This is obviously a multiple gap since the smaller ideals $I_i$ were already a multiple gap. By Theorem~\ref{Giarejigsaw} it is indeed a jigsaw.$\qed$\\

\begin{cor}
If $\{I_i : i\in n\}$ is a strongly countably separated dense $n$-gap in $\mathcal{P}(\omega)/fin$, then it is a jigsaw.
\end{cor}

Proof: By Corollary~\ref{containedindensejigsaw}, we have $I_i\subset J_i$ with $\{J_i: i\in n\}$ being a jigsaw. But this together with density implies that $\{I_i : i\in n\}$ is also a jigsaw. Namely, suppose $\{I_i|_a : i\in A\}$ is a multiple gap, and consider $B\subset A$. Then $\{J_i|_a : i\in A\}$ is also a multiple gap, hence there exists $b\subset a$ such that $\{J_i|_b : i\in B\}$ is multiple gap and $b\in J_i^\perp$ for $i\not\in B$. Then $b\in I_i^\perp$ for $i\not\in B$ since $I_i\subset J_i$. Also $\{I_i|_b : i\in B\}$ is a multiple gap. Indeed, if  $a_i\geq I_i|b$, then $a_i\geq J_i|_b$. For suppose $x\not\leq a_i$, with $x\in J_i$, $x\leq b$. Then by density, there exists $j$ and $z\in I_j$ with $z\subset x\setminus a_i$. Then $z\in J_i$ so we must have $j=i$. Hence $z\in I_i$ and $z\not\leq a_i$ a contradiction.$\qed$\\

\begin{cor}
If $\{I_i : i\in n\}$ is a strongly countably separated dense $n$-gap in $\mathcal{P}(\omega)/fin$, then there exists $a\subset\omega$ and a strongly countably separated jigsaw $\{J_i : i\in n\}$ made of Borel $F_{\sigma\delta}$ ideals such that $I_i|_a\subset J_i$ for every $i\in n$.\\
\end{cor}

\subsection{Example of a strongly countably separated clover}

Let $T=2^{<\omega}$ be the dyadic tree. We work in $\mathcal{P}(T)/fin$. Every element $x\in 2^\omega$ is viewed as a branch of $T$, in particular it is an infinite subset of $T$. Given $S\subset 2^\omega$, we call $I_S$ to the ideal of $\mathcal{P}(T)/fin$ generated by $S$. A Bernstein set is a subset $S\subset 2^\omega$ that uncountably meets every uncountable Borel subset of $2^\omega$. Bernstein sets can be constructed by transfinite induction by enumerating $2^\omega$ and its Borel subsets by the ordinal $\mathfrak c$. The classical Bernstein's result states that $2^\omega$ can be divided into two Bernstein sets, but the same kind of argument can provide a decomposition into any finite number of Bernstein sets.

\begin{prop}
Let $\{S_i : i\in n\}$ be pairwise disjoint Bernstein subsets of $2^\omega$. Then the ideals $\{I_{S_i} : i\in n\}$ constitute a clover $n$-gap which is strongly countably separated.
\end{prop}

Proof: It is clear that the ideals are mutually orthonal since $x\cap y =^\ast \emptyset$ for every branches $x\neq y$. If the ideals were not a multiple gap, then we would we have elements $c_i\geq I_{i}$ with $\bigcap_{i\in n}c_i =^\ast \emptyset$. Consider the sets $Z_i = \{x\in [T] : x\subset^\ast c_i\}\subset 2^\omega$ (remember that $[T]$ is naturally identified with $2^\omega$). The sets $Z_i$ are Borel and $S_i\subset Z_i$. Hence each $Z_i$ is co-countable and in particular $\bigcap_{i\in n}Z_i\neq\emptyset$. On the other hand $\bigcap_{i\in n}Z_i=\emptyset$ because $\bigcap_{i\in n}c_i =^\ast \emptyset$, a contradicition. The fact that the ideals are strongly countably separated is because if we put on $[T] = 2^\omega$ the Cantor set topology and we identify each $t=(t_0,\ldots,t_k)\in T$ with the branch $(t_0,t_1,\ldots,t_k,0,0,\ldots)$, then for any $S\subset[T]$ we have $I_S\subset\mathcal{I}_S$, so we can use Theorem~\ref{caracterizationcs}. In order to show now that the ideals are a clover, suppose that $d\subset T$ is such that $d\in I_{S_j}^\perp$ for some $j$. The set
$X=\{x\in [T] : x\cap d\neq^\ast\emptyset\}$ is a Borel set that is disjoint from $S_j$. Since $S_j$ is a Bernstein set, it follows that $X$ is countable. Therefore each ideal $I_{S_i}|_d$ is countably generated for $i\neq j$. Countably generated orthogonal ideals in $\mathcal{P}(\omega)/fin$ can always be separated.$\qed$\\

\section{Multiple gaps made of small ideals}\label{noHausdorff}

Given a cardinal $\theta$, we say that an ideal is $\theta$-generated if it is generated by a set of cardinality $\theta$. Given $k\in\omega$, a partialy ordered set (poset, for short) $\mathbb{P}$ is $\sigma$-$k$-linked if $\mathbb{P} = \bigcup_{m<\omega} P_m$ in such a way that for every $m$, every $k$ many elements of $P_m$ are compatible (that is, for every $\{p_i : i<k\}\subset P_m$ there exists $p\in \mathbb{P}$ such that $p\leq p_i$ for every $i\in k$). We write MA$_\theta$($\sigma$-$k$-linked) for Martin's axiom for cardinal $\theta$ and $\sigma$-$k$-linked posets: In a $\sigma$-$k$-linked poset, every family of $\theta$ many dense subsets has a generic filter. This is weaker than the usual MA$_\theta$, the same statement but for ccc partial orders. Let $\frak X$ be a family of subsets of $n$ which is upwards closed, that is: if $A\in\frak X$ and $B\supset A$, then $B\in\frak X$. We call $\frak X[k]$ to the family of all $A\subset n$ such that whenever $A=\bigcup_{j\in k}B_j$ then there exists $j\in k$ such that $B_j\in\mathfrak{X}$. One example is given by $\frak X = \{A\subset n : |A|\geq 2\}$ and $\frak X[k] = \{A\subset n : |A|> k\}$. All along this section, $\theta$ is a cardinal and $k\geq 1$ an integer. We notice that the results of this section hold in ZFC when we take $k=1$ and $\theta=\aleph_0$.

\begin{thm}[MA$_\theta$($\sigma$-$k$-linked)]\label{separationsmallideals}
Let $\{I_i : i\in n\}$ be $\theta$-generated ideals of $\mathcal{P}(\omega)/fin$ and $\mathfrak{X}$ an upwards closed family of subsets of $n$ such that $\bigcap_{i\in A}x_i =^\ast \emptyset$ whenever $A\in\frak X$ and $x_i\in I_i$ for every $i$. Then there exist $c_i\subset\omega$, $i\in n$ such that $I_i\leq c_i$ for $i\in n$, and such that $\bigcap_{i\in A} c_i = \emptyset$ for $A\in\frak X[k]$.
\end{thm}

Proof: We consider a poset $\mathbb{P}$ for application of MA$_\theta$($\sigma$-$k$-linked). The elements of $\mathbb{P}$ are of the form $p = (x_p^i, i\in n)$ where $x_p^i\in I_i$ for every $i\in n$ and $\bigcap_{i\in A}x_p^i = \emptyset$ for every $A\in\frak X[k]$. The order relation is given by $p\leq q$ if and only if $x_p^i \supset x_q^i$ for all $i\in n$.\\

\emph{Claim}: $\mathbb{P}$ is $\sigma$-$k$-linked. Proof of the claim: Given $\{u^{B} : B\in\mathfrak{X}\}$ finite subsets of $\omega$, let $Q$ be the set of all $p\in\mathbb{P}$ such that $u^{B} = \bigcap_{i\in B} x_p^i$ for every $B\in\frak X$. We prove that every $p[0],\ldots,p[k-1]\in Q$ are compatible. For this we show that if we define $r_i = \bigcup_{j\in k} x_{p[j]}^i$, then $r=(r_i)_{i\in n}\in\mathbb{P}$ and it satisfaies $r\leq p[j]$ for every $j\in k$. The second part is obvious, so we only have to check that $r\in\mathbb{P}$. So we fix $A\in\mathfrak{X}[k]$ and we prove that $\bigcap_{i\in A} r_i=\emptyset$. Suppose that we have $m\in \bigcap_{i\in A} r_i$. Call $B_j=\{i\in A : m\in x_{p[j]}^i\}$. Then $A=\bigcup_{j\in k}B_j$, so since $A\in\mathfrak{X}[k]$, there exists $j_0$ such that $B_{j_0}\in\mathfrak{X}$. For simplicity we suppose $j_0=0$. Keep in mind that, since $\mathfrak X$ is upwards closed, this implies that $B_0\cup\cdots\cup B_l\in\mathfrak{X}$ for every $l$. Then

\begin{eqnarray*}
m &\in& \bigcap_{i\in B_0}x_{p[0]}^i \cap \bigcap_{i\in B_1}x_{p[1]}^i \cap \cdots \cap \bigcap_{i\in B_{k-1}}x_{p[k-1]}^i \\
&=& u^{B_0} \cap \bigcap_{i\in B_1}x_{p[1]}^i \cap \cdots \cap \bigcap_{i\in B_{k-1}}x_{p[k-1]}^i\\
&=& \bigcap_{i\in B_0}x_{p[1]}^i \cap \bigcap_{i\in B_1}x_{p[1]}^i \cap \bigcap_{i\in B_2}x_{p[2]}^i \cap \cdots \cap \bigcap_{i\in B_{k-1}}x_{p[k-1]}^i\\
&=& u^{B_0\cup B_1} \cap \bigcap_{i\in B_2}x_{p[2]}^i \cap \bigcap_{i\in B_3}x_{p[3]}^i \cap \cdots \cap \bigcap_{i\in B_{k-1}}x_{p[k-1]}^i\\
&=& \cdots = u^{B_0\cup\cdots\cup B_{k-1}} = u^A = \emptyset
\end{eqnarray*}
This contradiction finishes the proof of the claim.\\

Now, for every $x=(x^i, i\in n)\in \prod_{i\in n}I_i$ , the set $$D_x = \{p\in\mathbb{P} : x^i\subset^\ast x_p^i,\ i\in n\}$$ is a dense subset of $\mathbb{P}$. Since the ideals are $\theta$-generated, we have cofinal subsets $S_i\subset I_i$ with $|S_i|=\theta$. By $MA_\theta$($\sigma$-$k$-linked), there is a filter $G$ in $\mathbb{P}$ which is generic for all the dense sets $D_x$, $x\in\prod_{i\in n}S_i$. This implies that $G$ is generic for all the dense set $D_x$, $x\in\prod_{i\in n}I_i$.  The sets $c_i=\bigcup\{x_p^i : p\in G\}$, for $i\in n$ are as desired.$\qed$\\

\begin{cor}[MA$_\theta$($\sigma$-$k$-linked)]\label{nogaps}
For $n>k$, there exist no $n$-gaps of $\theta$-generated ideals in $\mathcal{P}(\omega)/fin$.
\end{cor}

Proof: Let $\{I_i : i\in n\}$ be mutually orthogonal ideals. Then the hypotheses of Theorem~\ref{separationsmallideals} hold for $\frak X = \{ A\subset n : |A|\geq 2\}$, hence there exist elements $c_i\geq I_i$ with $\bigcap_{i\in A}c_i = \emptyset$ for all $A\in \frak X[k] = \{A\subset n : |A|>k\}$.$\qed$\\

Next corollary is just the topological interpretation of Theorem~\ref{separationsmallideals}. An open subset of a topological space $X$ is called $F_\theta$ if it is the union of $\theta$ many closed sets. We point out that, by Lemma~\ref{Xseparation} the conclusion of Corollary~\ref{separationsmallopen} is equivalent to say that $\bigcap_{i\in A}\overline{U_i}=\emptyset$ for all $A\in\frak X[k]$.

\begin{cor}[MA$_\theta$($\sigma$-$k$-linked)]\label{separationsmallopen}
Let $\{U_i : i\in n\}$ be open $F_\theta$ subsets of $\omega ^\ast$, and $\mathfrak{X}$ an upwards closed family of subsets of $n$ such that $\bigcap_{i\in A}U_i = \emptyset$ for all $A\in \frak X$. Then there exist clopen sets $C_i\supset U_i$ such that $\bigcap_{i\in A}C_i=\emptyset$ for all $A\in\frak X[k]$.\\
\end{cor}

We present a result for the Banach lattice $C(\omega^\ast)$ of real-valued continuous functions on $\omega^\ast$ (that is, the Banach lattice $\ell_\infty/c_0$), which is a consequence of Theorem~\ref{separationsmallideals}. We use the following notation: Given $\mathcal{F}\subset C(\omega^\ast)$ and $g\in C(\omega^\ast)$, $g\leq \mathcal{F}$ means that $g\leq f$ for all $f\in\mathcal{F}$.\\

\begin{thm}[MA$_\theta$($\sigma$-$k$-linked)]\label{separationsmallfunctions}
Let $\{\mathcal{F}_i\geq 0 : i<n\}$, be $n$ subsets of $C(\omega ^\ast)$ of cardinality $\theta$ such that $\sum_{i\in n}f_i\leq 1$ whenever $f_i\in\mathcal{F}_i$. Then, there exist functions $g_i\in C(\omega ^\ast)$, such that $g_i\geq\mathcal F_i$ and $\sum_{i\in n}g_i\leq k$.
\end{thm}

Proof: As a first step, we prove the following claim:

Claim 1: Given $\varepsilon>0$ there exist functions $g_i\geq \mathcal{F}_i$ such that $\sum_{i\in n}g_i\leq k+\varepsilon$.
Proof of the claim: Consider $m$ such that $n/2^m<\varepsilon$. Let $\Delta = \{j/2^m : j=0,\ldots,2^m-1\}$. For every $\delta\in \Delta$ and $i\in n$, let $$U_\delta^{i} = \{x\in \omega ^\ast : \exists f\in \mathcal{F}_i \ f(x)>\delta\}$$

These are $F_\theta$ open subsets of $\omega ^\ast$. Given $\delta_\ast = (\delta_0,\ldots,\delta_{n-1})\in \Delta^n$, we consider $\frak X_{\delta_\ast} = \{A\subset n : \sum_{i\in A}\delta_i\geq 1\}$. Notice that
$\frak X_{\delta_\ast}[k] \supset \{A\subset n : \sum_{i\in A}\delta_i\geq k\}$. By Corollary~\ref{separationsmallopen} we obtain clopen sets $C_{\delta_\ast}^i \supset U_{\delta_i}^i$ such that $\bigcap_{i\in A}C_{\delta_\ast}^i = \emptyset$ whenever $A\in\frak X_{\delta_\ast}[k]$.
For $\delta\in\Delta$, $i\in n$, set $D_\delta^i =\bigcap\{C_{\delta_\ast}^i : \delta_i = \delta\}$. Observe that $D_\delta^i \supset U_\delta^i$. Define then
$$g_i(x) = 2^{-m} + \max\{\delta\in\Delta : x\in D_\delta^i\}$$

We define that maximum to be 0 if $x\not\in\bigcup_{\delta\in\Delta}D_\delta^i$. We check the desired properties:
\begin{itemize}
\item The functions $g_i$ are continuous because the sets $D_\delta^i$ are finitely many clopen sets.
\item $g_i\geq \mathcal{F}_i$. Let $f\in\mathcal{F}_i$ and $x\in\omega ^\ast$ with $f(x)>0$. Pick $\delta\in\Delta$ such that $\delta < f(x)\leq \delta + 2^{-m}$. Then $x\in U_\delta^i\subset D_\delta^i$, hence $g_i(x) \geq 2^{-m} + \delta \geq f(x)$.
\item $\sum_{i\in n} g_i\leq k+\varepsilon$. Pick $x\in\omega ^\ast$. Say that $g_i(x) = 2^{-m} + \delta_i$, so that $x\in D_{\delta_i}^i$ for $i\in A =\{j : \delta_j\neq 0\}$. Consider $\delta_\ast = (\delta_0,\ldots,\delta_{n-1})$. We have that $$x\in\bigcap_{i\in A} D_{\delta_i}^i \subset \bigcap_{i\in A} C_{\delta_\ast}^i,$$ so in particular that intersection is nonempty, so $A\not\in \frak X_{\delta_\ast}[k]$, hence $\sum_{i\in n}\delta_i = \sum_{i\in A}\delta_i<k$. Therefore $\sum_{i\in n} g_i(x) < n2^{-m} + k < \varepsilon + k$.
\end{itemize}

This finishes the proof of Claim 1.\\

Claim 1 allows to find functions $g_i^m$ for $i\in n$ and $m\in\omega $ such that $g_i^m\geq \mathcal{F}_i$ and $\sum_{i\in n}g_i^m \leq k+2^{-m}$. We can take these functions so that $g_i^0\geq g_i^1\geq g_i^2\geq\cdots$. Moreover, the hypotheses of the theorem imply that $0\leq\mathcal{F}_i\leq 1$ for each $i$, so we can suppose without loss of generality that $g^m_i:\omega^\ast\To [0,1]$ (we could change $g^m_i$ by $\min(g^m_i,1)$). To finish the proof, we find functions $g_i\in C(\omega^\ast)$ such that $g_i\geq \mathcal{F}_i$ and $g_i\leq g_i^m$ for all $m$. This is a consequence of the fact that, under MA$_\theta$($\sigma$-centered) there are no $(\aleph_0,\theta)$-gaps in $\mathcal{P}(\omega)/fin$, since gaps in $\mathcal{P}(\omega)/fin$ are related to gaps in $C(\omega^\ast)$. Anyway, we provide a direct proof. We show that if $\mathcal{F}\geq 0$ has cardinality $\theta$ and $g^1\geq g^2\geq\cdots \geq \mathcal{F}$, then there exists $g$ such that $g^1\geq g^2\geq\cdots\geq g \geq \mathcal{F}$. We consider continuous functions $\hat{g}_m:\beta\omega\To [0,1]$ such that $\hat{g}_m|_{\omega^\ast} = g^m$. Notice that for $\phi,\psi\in C(\beta\omega)$ we have $\phi|_{\omega^\ast}\leq \psi|_{\omega^\ast}$ if and only if $\phi(k)\leq\psi(k)$ for all but finitely many $k\in\omega$. Consider a poset $\mathbb{P}$ whose elements are pairs of the form $p=(h^p,s^p)$ where $h^p:\beta\omega\To [0,1]$ is continuous and $s^p\in\omega^{<\omega}$ with the following properties: \begin{itemize} \item $h^p|_{\omega^\ast}\leq g^m$ for every $m<\omega$, \item $h^p(k) \leq \hat{g}_m(k)$ for every $m<length(s^p)$ and every $k>s^p_m$.
\end{itemize}
The ordered relation is that $p\leq q$ if $h_q\leq h_p$ (as functions on $\beta\omega$) and $s^q\leq s^p$ (as elements of the tree $\omega^{<\omega}$). This partial order is $\sigma$-$k$-linked (even $\sigma$-centered) since a finite number of conditions with a fixed $s^p$ are compatible. We have a filter $G$ that is generic for each of the dense sets $D_m = \{p : length(s_p) >m\}$ for $m\in\omega$ and $D_f = \{p : h_p|_{\omega^\ast}\geq f\}$ for $f\in\mathcal{F}$. Let $g(k) = \sup_{p\in G}h_p(k)$ for $k\in\omega$ and $g:\beta\omega\To [0,1]$ a continuous extension. The function $g|_{\omega^\ast}$ is the desired one.$\qed$\\

We notice that if there exists a $k$-gap of $\theta$-generated ideals in $\mathcal{P}(\omega)/fin$, then the constant $k$ of Theorem~\ref{separationsmallfunctions} cannot be improved. Namely, consider $\{I_i : i\in k\}$ a multiple gap, and $\mathcal{F}_i\subset C(\omega^\ast)$ the set of characteristic functions of $\theta$ many clopens whose union is $U(I_i)$. Suppose that we have $g_i\geq \mathcal{F}_i$. Pick $x\in\bigcap_{i\in k}\overline{U(I_i)}$. Then $g_i(x)\geq 1$ for every $i$, hence $\sum_{i\in k}g_i(x) \geq k$.\\

\begin{thm}
Let $M$ be a model of ZFC. Let $\aleph_0<\theta_1< \theta_2 < \cdots$ be cardinals in $M$ and $1<n_1\leq n_2 \leq \cdots$ be integers. Then, there exists a ccc generic extension $M(G)$ of $M$ such that for every $j\in\omega$ the following hold in $M(G)$:
\begin{enumerate}
\item MA$_{\theta_j}$($\sigma$-$n_j$-linked), so there are no $\theta_j$-generated $n$-gaps for $n>n_j;$
\item There exists an $n_j$-gap of $\theta_j$-generated ideals in $\mathcal{P}(\omega)/fin$.
\end{enumerate}
\end{thm}

Proof: Let $M_1$ be the generic extension of $M$ obtained by adding $\theta_\omega$ many Cohen reals, where $\theta_\omega$ is any cardinal larger than all cardinals $\theta_j$. For every $j$, we can consider a set of Cohen reals $X_j \subset n_j^\omega$ of cardinality $\theta_j$. We identify the set $n_j^\omega$ with the set of branches of the $n_j$-adic tree $n_j^{<\omega}$. A set $X$ of Cohen reals is always a Luzin set, that is, every nowhere dense subset of $X$ is countable. This implies in particular that each of the sets $X_j$ has property $A(\aleph_1,n_j)$ described as follows:\\

Definition: Let $\theta$ be a cardinal, $n\in\omega$, and $X\subset n^\omega$. We say that $X$ has property $A(\theta,n)$ if for every $Y\subset X$ with $|Y|=\theta$, there exists $t\in n^{<\omega}$ and elements $\{y^i : i\in n\}\subset Y$ such that $t^\frown i \in y^i$ for every $i\in n$.\\

In $M_1$ pick a cardinal $\kappa\geq\theta_\omega$ such that $\kappa^{\theta_\omega}=\kappa.$ Let $M_2$ be the generic extension of $M_1$ obtained by a a natural finite support iteration $(\mathbb{P}_\alpha, \mathbb{Q}_\alpha: \alpha<\kappa)$ of posets whose limit $\mathbb{P}_\kappa$ force MA$_{\theta_j}$($\sigma$-$n_j$-linked) for every $j$. Thus, for each of the iterands $\mathbb{Q}_\alpha$ of the iteration there is $j$ so that $\mathbb{P}_\alpha$ forces both that $|\mathbb{Q}_\alpha| = \theta_j$ and that $\mathbb{Q}_\alpha$ is $\sigma$-$n_j$-linked. The model $M_2$ is the one we are looking for. It remains to find a $n_j$-gap of $\theta_j$-generated ideals in the model $M_2$. For every $j$ we fix a regular cardinal $\lambda_j$ such that $\theta_{j-1}<\lambda_j \leq \theta_j$ (assume $\theta_0 = \aleph_0$).\\

Definition. Given a regular cardinal $\lambda$ and $n\in\omega$, we say that a poset $\mathbb{P}$ has property $K(\lambda,n)$ if every subset of $\mathbb{P}$ of cardinality $\lambda$ has a further subset of cardinality $\lambda$ which is $n$-linked. \\

\noindent It is well-known and easily seen that this Knaster-type chain conditions are preserved under finite support iterations.
If a poset $\mathbb{Q}$ satisfies either $|\mathbb{Q}|<\lambda$ or it is $\sigma$-$n'$-linked for $n'\geq n$, then $\mathbb{Q}$ has property $K(\lambda,n)$. It follows that the poset that forces the generic extension from $M_1$ to $M_2$ has property $K(\lambda_j,n_j)$ for every $j$.\\

Claim 1: Let $n\in\omega$ and $\lambda$ be a regular cardinal. Let $X\subset n^\omega$ have propery $A(\lambda,n)$. Let $\mathbb{P}$ be a poset that satisfies $K(\lambda,n)$. Then the set $X$ still has property $A(\lambda,n)$ in the generic extension forced by $\mathbb{P}$. \\

\noindent Proof of the claim: Suppose not. Then
$$M(G)\models \exists Y\subset X\ |Y|=\lambda\text{ and } \forall t\in n^{<\omega}\ \exists \iota(t)\in n\ \forall y\in Y\ t^\frown \iota(t)\not\in y$$
Let $\dot{Y}$ and $i$ be names for $Y$ and the function $\iota:n^{<\omega}\To n$ respectively. Let $\tilde{Y} = \{x\in X : \exists p\in\mathbb{P} : p\Vdash x\in\dot{Y}\}$. Then, $\tilde{Y}$ is a set in the ground model with $Y\subset \tilde{Y}\subset X$. Since the forcing is $K(\lambda,n)$, in particular $\lambda$-cc, cardinal $\lambda$ is preserved, hence $|\tilde{Y}|\geq \lambda$. For every $x\in \tilde{Y}$ there exists $p(x)\in \mathbb{P}$ such that
\begin{eqnarray*}
p(x) &\Vdash& \dot{Y}\subset X,\ |\dot{Y}| = \lambda \text{ and } x\in\dot{Y}\\
p(x) &\Vdash& \forall t\in T\  i(t)\in n\\
p(x) &\Vdash& \forall t\in T\ t^\frown i(t)\not\in x
\end{eqnarray*}

Since $\mathbb{P}$ has property $K(\lambda,n)$ and $|\tilde{Y}|\geq\lambda$ is a regular cardinal, there exists $Z\subset\tilde{Y}$ with $|Z| = \lambda$ such that $\{p(x) : x\in Z\}$ is $n$-linked. Since $X$ has property $A(\lambda,n)$ there exist $t\in n^{<\omega}$ and $\{x^j : j\in n\}\subset X$ such that $t^\frown j\in x^j$ for all $j\in n$. Since $Z$ is $n$-linked, we can pick a condition $p\leq p(x_j)$ for $j\in n$. Then $p\Vdash t^\frown i(t)\not\in x^j$, hence $p\Vdash i(t)\neq j$ for every $j\in n$ which contradicts the fact that $p\Vdash i(t)\in n$. This finishes the proof of Claim 1.\\

From Claim 1 and the remarks before it, we get that $M_2\models X_j$ has property $A(\lambda_j,n_j)$ for every $j$. It remains to show the following:\\

Claim 2: Let $n\in\omega$ and $\lambda$ a regular cardinal. Let $X\subset n^\omega$ with property $A(\lambda,n)$. Like in Section~\ref{analyticjigsaw}, for every $x\in X$ and $i\in n$ consider $a_x^i = \{s\in T : s^\frown i\in x\}$. Let $J_i$ be the ideal generated by $\{a_x^i : x\in X\}$. Then the ideals $\{J_i : i\in n\}$ form a $n$-gap that is not weakly countably separated (indeed not even weakly $(<\lambda)$-separated). Proof of the claim: We just follow the same argument as for Theorem~\ref{analyticjigsawexample}. Let us denote by $T(m)$ the set of all elements of $t\in T =n^\omega$ of length less than $m$. Suppose that the ideals were weakly countably separated, so that we have elements $c_i^k\subset \omega$, $i\in n$, $k\in \omega$, such that $\bigcap_{i\in n}c_i^k = \emptyset$ for each $k$, and whenever we pick $b_i\in\mathcal{J}_i$ there exists $k$ with $b_i\subset^\ast c_i^k$. In particular, for every $x\in X$ and every $i\in n$ there exist $k(x),m(x)\in \omega$ such that $a_x^i \setminus T(m(x)) \subset c_i^{k(x)}$ for every $i\in n$. For every $m,k\in \omega$ and every $s\in n^{m+2}$, let $$X(m,k,s) = \{x\in X  : x|_{m+2} = s\text{ and }\forall i\in n\ a_x^i \setminus T(m) \subset c_i^{k}\}.$$ There exist $m_0,k_0,s_0$ such that $|X(m_0,k_0,s_0)|\geq \lambda$, and by property $A(\lambda,n)$ we can find $t\in T$ (necessarily of length greater than $m_0$) and $\{x^i : i\in n\}\subset X(m_0,k_0,s_0)$ with $t^\frown i \in x^i$. This implies that $t\in\bigcap_{i\in n}a_{x^i}^i \setminus T(m_0)\subset\bigcap_{i\in n}c_i^{k_0}=\emptyset$, a contradiction.$\qed$\\

\section{Dense multiple gaps from completely separable almost disjoint families}\label{completelyseparablesection}

An almost disjoint family is a family $\mathcal{A}$ of infinite subsets of $\omega$ such that $a\cap b =^\ast \emptyset$ for all different $a,b\in\mathcal{A}$.

\begin{defn}
An almost disjoint family $\mathcal{A}$ is called completely separable if for every subset $x\subset\omega$ that is not in the ideal generated by $\mathcal{A}$ we have that $|\{a\in\mathcal{A}: a\subset^\ast x\}|=\mathfrak c$.\\
\end{defn}

The concept traces back to \cite{Hechler}. Completely separable almost disjoint families exist under the assumption that $\mathfrak a = \mathfrak c$ \cite{Rothberger} or $\mathfrak s=\omega_1$ \cite[Corollary 2.8]{BalDocSim}. A recent work of Shelah \cite{Shelahmadsane} shows their existence under other assumptions, in particular if $\mathfrak s<\mathfrak a$ or $\mathfrak c<\aleph_\omega$. The consistency of the non existence of such families is unknown.\\

\begin{thm}
Assume a completely separable almost disjoint family exists. Let $n\in\omega$ and $\mathfrak X$ be a family of proper nonempty subsets of $n$. Then, there exists a dense $n$-gap $\{I_i : i\in n\}$ in $\mathcal{P}(\omega)/fin$ that is a $B$-jigsaw for $B\in\mathfrak X$ but a $B$-clover for $B\not\in\mathfrak X$ $(|B|\geq 2)$.
\end{thm}

Proof: We proceed by induction on $n$, so we fix $n$ and we suppose the theorem proven for all $m<n$. For convenience we will assume that $\mathfrak{X}$ includes all singletons. For the case $n=1$ we just need an ideal and for $n=2$ we just need a 2-gap. Let $\mathcal{A}$ be a completely separable almost disjoint family. We decompose our family into disjoint subfamilies in the form $\mathcal{A} = \bigcup_{B\in\mathfrak X}\mathcal{A}_B$ in such a way that if $x\subset\omega$ is not in the ideal generated by $\mathcal{A}$, then $|\{a\in\mathcal{A}_B : a\subset^\ast x\}|=\mathfrak c$ for all $B\in\mathfrak X$. This decomposition can be made as follows: let $\{x_\alpha : \alpha<\mathfrak c\}$ be an enumeration of all infinite subsets of $\omega$ that are not in the ideal generated by $\mathcal A$ in which each set appears repeated $\mathfrak c$ many times. Inductively on $\alpha$, we will choose elements $\{a_B^\alpha : \alpha< \mathfrak c, B\in\mathfrak{X}\}\subset\mathcal{A}$. At step $\alpha$, since $|\{a\in\mathcal{A} : a\subset^\ast x\}|=\mathfrak c$ and we have chosen only less than $\mathfrak c$ elements in the previous steps, we can pick $a^\alpha_B\subset^\ast x_\alpha$ and different from all previous choices. At the end, it is enough to declare $a^\alpha_B\in\mathcal{A}_B$ for all $\alpha<\mathfrak c$ and all $B\in\mathfrak X$ and we will have a decomposition as required.\\

Now we define the ideals $\{I_i : i\in n\}$. For $B\in\mathfrak X$, $2\leq |B|<n$ and $x\in\mathcal{A}_B$, we will apply the inductive hypothesis to find a dense multiple gap $\{I^x_i : i\in B\}$ in the Boolean algebra $\mathcal{P}(x)/fin$ that is a $C$-jigsaw for $C\subset B$, $C\in\mathfrak X$ but $C$-clover for $C\subset B$, $C\not\in\mathfrak X$. For every $i\in n$ the ideal $I_i$ is defined as the ideal generated by
$$\mathcal{A}_{\{i\}}\cup\bigcup\left\{I_i^x : x\in\bigcup\{\mathcal{A}_B : B\in\mathfrak X, i\in B, 2\leq |B| < n\}\right\}.$$

Let us prove that this is a dense multiple gap. Density is clear. The ideals are mutually orthogonal, because if $u$ and $v$ belong to the set of generators of $I_i$ and $I_j$ respectively ($i\neq j$), then:  either $u\subset x$ and $v\subset y$ for different elements $x,y\in\mathcal{A}$ (and $\mathcal{A}$ is almost disjoint); or there exists $x\in\mathcal{A}$ such that $u\in I_i^x$ and $v\in I_j^x$. Finally, suppose that we have elements $c_i$ with $c_i\geq I_i$. We prove inductively on $k<n$ that $\bigcap_{i\in k} c_i$ is not in the ideal generated by $\mathcal{A}$, in particular it is infinite. So suppose $x=\bigcap_{i\in k} c_i$ is not in the ideal generated by $\mathcal{A}$, hence there are $\mathfrak c$ many $a\in \mathcal{A}_{\{k\}}$ such that $a\subset^\ast x$. On the other hand, $a\subset^\ast c_{k}$ for all $a\in\mathcal{A}_{\{k\}}$, so there are $\mathfrak c$ many $a\in\mathcal{A}_{\{k\}}$ with $a\subset^\ast c_k\cap x = \bigcap_{i\in k+1}c_i$. This proves that the ideals form a multiple gap.\\

Let us check now that the ideals form a $B$-jigsaw for $B\in\mathfrak X$. So assume we have $a\subset \omega$ and $A\supset B$ such that $\{I_i|_a : i\in A\}$ form a multiple gap. There are two possibilities: either $a$ is in the ideal generated by $\mathcal{A}$ or it is not. If it is not, then $a$ contains an element $b\in \mathcal{A}_B$. Then $b\in I_i^\perp$ for $i\not\in B$ while $\{I_i|_b : i\in B\} = \{I^b_i : i\in B\}$ is a multiple gap, so we are done. The other case is that $a$ is in the ideal generated by the family $\mathcal{A}$. Let $F = \{u\in \mathcal{A} : u\cap a \neq^\ast \emptyset\}$. This set is finite and $a = \bigcup_{u\in F}  a\cap u$. Since we assume that $\{I_i|_a : i\in A\}$ is a multiple gap, there must exist $u\in F$ with $\{I_i|_{a\cap u} : i\in A\} = \{I_i^u|_{a\cap u} : i\in A\}$ is a multiple gap. For that $u$, we must have $u\in\mathcal{A}_{C}$ with $C\supset A$. Since we assumed that the ideals $\{I_i^u : i\in A\}$ are a $B$-jigsaw, we find $b\subset a\cap u$ such that $b\in I_i^\perp$ for $i\in A\setminus B$ but $\{I_i|_b : i\in B\}$ is a multiple gap. This finishes the proof that the ideals are a $B$-jigsaw for $B\in\mathfrak X$.\\

We prove that the ideals are a $B$-clover for $B\not\in\mathfrak X$. We suppose for contradiction that there exists an infinite set $a\subset \omega$ such that $\{I_i|_a : i\in B\}$ is a multiple gap but $a\in I_i^\perp$ for $i\not\in B$. Again we distinguish two cases, according whether $a$ in the ideal generated by $\mathcal A$ or not. If $a$ is not in the ideal, then it contains $\mathfrak c$ many elements from $\mathcal{A}_{\{i\}}$ for every $i\in n$, and that contradicts that $a\in I_i^\perp$ for $i\not\in B$. So we suppose that $a$ is in the ideal generated by $\mathcal{A}$ and again we can consider the finite set $F=\{u\in\mathcal{A} : u\cap a \neq^\ast 0\}$, so that $a=\bigcup_{u\in F}a\cap u$. Since $\{I_i|_a : i\in B\}$ is a multiple gap, again we must have that for some $u\in F$, $\{I_i|_{a\cap u}:i\in B\}$ is a multiple gap, so $u\in \mathcal{A}_C$ with $C\supset B$. Notice that $C\neq B$ since $B\not\in \mathfrak X$. We have that $\{I_i|_{a\cap u}:i\in B\} = \{I_i^u|_{a\cap u} : i\in B\}$ is a multiple gap but $a\cap u\in (I_i^u)^{\perp}$ for $i\in C\setminus B$. This contradicts the hypothesis that the ideals $\{I_i^u : i\in C\}$ are a $B$-clover.$\qed$\\

\section{Dense jigsaws and injective Banach spaces}\label{injectiveBanach}

In this section we consider an application to the theory of Banach spaces. We refer to \cite{AlbKal,Banach} for the basic facts and terminology on the subject.\\

The existence of a gap in a Boolean algebra is equivalent to the failure of completeness. Namely, in a complete Boolean algebra orthogonal ideals can always be separated by their suprema, and conversely if a set $S$ lacks a supremum, the ideal generated by $S$ forms a gap together with its orthogonal ideal. At the same time, completeness is equivalent to injectivity (cf. for instance \cite{Walker}), meaning that $\mathcal{B}$ is complete if and only if for every superalgebra $\mathcal{B}'\supset\mathcal{B}$ there is a Boolean projection from $\mathcal{B}'$ onto $\mathcal{B}$. If we look at the Banach space of continuous functions $C(St(\mathcal{B}))$ we have in a similar spirit that $\mathcal{B}$ is complete if and only if $C(St(\mathcal{B}))$ is 1-injective, meaning that for every superspace $X\supset C(St(\mathcal{B}))$ there exists a projection $T: X\To C(St(\mathcal{B}))$ of norm~1 (cf. for instance \cite{AlbKal}). The situation becomes more complicated if we deal with injectivity instead of 1-injectivity of Banach spaces. The definition of injectivity is the same as 1-injectivity but we require the projection $T$ to be just bounded operator, not necessarily of norm 1. A characterization of injective Banach spaces is an important open problem. In this section we prove that while the existence of a gap is necessary and sufficient condition for the failure of 1-injectivity, the existence of arbitrarily large dense jigsaws is sufficient for the failure of injectivity.

\begin{thm}\label{noninjectivityboolean}
If $\mathcal{B}$ contains a dense $n$-jigsaw for every $n$, then $C(St(\mathcal{B}))$ is uncomplemented in a superspace $X$ of density character $|\mathcal{B}|$.
\end{thm}

This will be obtained as an application of a result by Ditor~\cite{Ditor}. Theorem~\ref{noninjectivityboolean} is a corollary of Theorem~\ref{noninjectivitygeneral}, whose proof is given at the end of this section, after some lemmas. For $\mathcal{B}=\mathcal{P}(\omega)/fin$, the assumption of Theorem~\ref{noninjectivityboolean} above holds by Theorem~\ref{Giarejigsaw} (when $\bigcup_{i\in n}G_i = K$). In this way, we have a proof of a result of Amir~\cite{Amir} that $C(\omega^\ast)$ is uncomplemented in a superspace of density character $\mathfrak c$. In relation with this, we mention the open problem whether $C(\omega^\ast)$ is uncomplemented in a superspace $X$ with $dens(X/C(\omega^\ast))=\aleph_1$. Theorem~\ref{noninjectivityboolean} can be stated more generally, removing the zero-dimensionality assumption. For this we define:

\begin{defn}
Let $K$ be a compact space and $\{F_i : i\in n\}$ closed subspaces. We say that they form a jigsaw of closed sets if for every $A\supset B$ nonempty proper subsets of $n$, the set $\bigcap_{i\in B}F_i\setminus\bigcap_{i\in A}F_i$ is dense in $\bigcap_{i\in B}F_i$. Such a jigsaw is called dense if $\bigcup_{i\in n}F_i = K$.
\end{defn}

A finite family of ideals $\{I_i : i\in n\}$ of $\mathcal{B}$ is a dense jigsaw if and only if the sets $\{\overline{U(I_i)}\}$ form a dense jigsaw of closed sets in $St(\mathcal{B})$. The more general statement is:

\begin{thm}\label{noninjectivitygeneral}
Let $K$ be a compact space. Suppose that for every $n$ there exists a dense jigsaw of $n$ many closed sets in $K$. Then $C(K)$ is uncomplemented in a superspace of the same density character as $C(K)$.
\end{thm}

Let $\{F_i : i\in n\}$ be a dense jigsaw of closed sets. We consider the compact space $L = \bigcup_{i\in n}F_i\times\{i\}$, the disjoint sum of the compact spaces $F_i$, and $\phi:L\To K$ the continuous surjection given by $\phi(x,i) = x$. The composition operator $\phi^0:C(K)\To C(L)$, given by $\phi^0(f) = f\circ\phi$, is an isometric embedding of Banach spaces. Let $H(L)$ be the hyperspace of $L$ consisting of all nonempty closed
subsets of $L$ endowed with the Vietoris topology. We define a derivation process on subsets of $K$, the derived set $D(X)$ of $X\subset K$ being defined as the set of all $x\in X$ for which there exist two disjoint closed nonempty sets $R,S\subset \phi^{-1}(x)$ such that both $R$ and $S$ belong to the closure in $H(L)$ of the set $\{\phi^{-1}(y) : y\in X\}$. As usual, the iterated derived sets are $D^{(n+1)}(X) =
D(D^{(n)}(X))$. This derivation procedure is defined by Ditor \cite[Definition 5.3]{Ditor}, our $D^{(n)}(K)$ would be in his notation $\Delta_\phi^{(n)}(2,2,\ldots,2)$. He proves the following~\cite[Corollary 5.4]{Ditor}:

\begin{lem}
If $D^{(k)}(K)\neq\emptyset$ and $T:C(L)\To C(K)$ is an operator with $T\phi^0 = 1_{C(K)}$, then $\|T\|\geq 1+k$.
\end{lem}

By $1_{C(K)}$ we denote the identity map on $C(K)$. In our particular situation, we have:

\begin{lem}
$D^{(k)}(K) = \left\{x\in K : |\{i\in n : x\in F_i\}|\geq 2^k\right\}$
\end{lem}

Proof: We prove the lemma by induction. It is obvious for $k=0$ ($D^{(0)}(K)=K$), so we suppose it holds for $k$ and prove it for $k+1$. Suppose $|\{i\in n : x\in F_i\}|\geq 2^{k+1}$. Put $\{i\in n : x\in F_i\} = A \cup B$ as the union of two disjoint sets with $|A|,|B|\geq 2^k$. Notice that $\phi^{-1}(x) = \{(x,i) : i\in A \cup B\}$. Because the sets $F_i$ form a jigsaw, $x$ is in the closure of the sets $V=\bigcap_{i\in A}F_i\setminus\bigcup_{i\not\in A}F_i$ and $W=\bigcap_{i\in B}F_i\setminus\bigcup_{i\not\in B}F_i$, which are both contained in $D^{(k)}(K)$ by the inductive hypothesis. Consider nets $x_\alpha\subset V$ and $y_\alpha\subset W$ that converge to $x$ and with $\phi^{-1}(x_\alpha)$ and $\phi^{-1}(y_\alpha)$ being convergent nets in the Vietoris topology to certain sets $R$ and $S$ respectively. Notice that $R\subset\{(x,i) : i\in A\}$ and $S\subset\{(x,i) : i\in B\}$. Hence this shows that $x\in D^{(k+1)}(K)$. For the converse, observe that for every $y\in D^{(k)}(K)$, $|\phi^{-1}(y)| = |\{(y,i) : y\in F_i\}|\geq 2^k$, and that for any $R$ in the closure of $\{\phi^{-1}(y) : y\in D^{(k)}(K)\}$ in the Vietoris topology we have $|R|\geq 2^k$. Therefore, if $|\phi^{-1}(x)|<2^{k+1}$, $x\not\in D^{(k+1)}(K)$.$\qed$\\

Putting together the two previous lemmas, we get:

\begin{lem}\label{finallemma}
If we have a jigsaw of $2^m$ many closed sets in $K$, and $L$ and $\phi$ are as above, then for any $T:C(L)\To C(K)$ with $T\phi^0=1_{C(K)}$ we have that $\|T\|\geq 1+m$.
\end{lem}

Proof of Theorem~\ref{noninjectivitygeneral}: By Lemma~\ref{finallemma}, for every $m$, we would have a continuous surjection $\phi_m:L_m\To K$ such that $L_m$ has the same weight as $K$, and for any $T:C(L_m)\To C(K)$ with $T\phi_m^0=1_{C(K)}$ we have that $\|T\|\geq 1+m$. Consider $$L = \left\{ (x_m)_{m<\omega}\in \prod_{m<\omega}L_m : \forall n,m<\omega\ \phi_m(x_m) = \phi_n(x_n)\right\}$$

The weight of $L$ equals the weight of $K$, and we have a continuous surjection $\phi:L\To K$ given by $\phi((x_m)_{m<\omega}) = \phi_n(x_n)$ for no matter what choice of a fixed $n$. $Y = C(L)$ will be the superspace that we look for, and we consider $C(K)\subset C(L)$ through the embedding $\phi^0:C(K)\To C(L)$. For every $m$ we can express $\phi = \phi_m \pi_m$ as the composition of $m$-th coordinate projection $\pi_m:L\To L_m$ with $\phi_m$, $\phi:L\To L_m\To K$. Hence we can write $C(K)\subset C(L_m) \subset C(L)$. Since any projection $T:C(L_m)\To C(K)$ must have norm at least $1+m$, no projection $T:C(L)\To C(K)$ exists.$\qed$\\

\end{document}